\documentclass{amsart}
\usepackage[nobysame]{amsrefs}
\usepackage{amssymb,mathrsfs}
\usepackage[all]{xy}
\newdir{ >}{{}*!/-5pt/@{>}} 
\SelectTips{cm}{} 
\usepackage{amssymb}
\usepackage{longtable}
\usepackage{array}
\usepackage{lscape}
\usepackage{tikz}

\allowdisplaybreaks 

\title[The Cohomology of the mod 2 Steenrod Algebra]
      {The Cohomology of the mod 2 Steenrod Algebra \\
	DOI:10.11582/2021.00077}
\author{Robert R. Bruner}
\address{Department of Mathematics, Wayne State University, USA}
\email{robert.bruner@wayne.edu}
\author{John Rognes}
\address{Department of Mathematics, University of Oslo, Norway}
\email{rognes@math.uio.no}

\date{February 2022}

\newtheorem{theorem}{Theorem}[section]
\newtheorem{proposition}[theorem]{Proposition}

\theoremstyle{definition}

\theoremstyle{remark}
\newtheorem{remark}[theorem]{Remark}

\numberwithin{equation}{section}

\DeclareMathOperator{\Ext}{Ext}
\DeclareMathOperator{\Hom}{Hom}
\DeclareMathOperator{\End}{End}

\newcommand{\bF}{\mathbb{F}}

\newcommand{\lla}{\longleftarrow}
\newcommand{\lra}{\longrightarrow}
\newcommand{\Lra}{\Longrightarrow}

\renewcommand{\:}{\colon}

\renewcommand{\>}{\rangle}
\newcommand{\cA}{\mathcal{A}}
\newcommand{\fm}{\mathfrak{m}}
\renewcommand{\Im}{\mathrm{Im}}
\newcommand{\Ker}{\mathrm{Ker}}

\begin{document}

\begin{abstract}
A minimal resolution of the mod 2 Steenrod algebra in the range $0
\leq s \leq 128$, $0 \leq t \leq 200$, together with chain maps for
each cocycle in that range and for the squaring operation $Sq^0$
in the cohomology of the Steenrod algebra.  
\end{abstract}

\maketitle

\noindent
This article describes the archived dataset~\cite{Ext-A-200},
available for download at the NIRD Research Data Archive {\tt
https://archive.sigma2.no}.  Please refer to the dataset and this article
by their digital object identifier {\tt DOI:10.11582/2021.00077}.

\setcounter{tocdepth}{1} 	
\tableofcontents

\section{Introduction}

Let $\cA$ denote the classical mod 2 Steenrod algebra.
This archive contains
\begin{enumerate}
\item
a minimal free resolution of $\bF_2$ over $\cA$ in internal degrees $t
\leq 200$ and cohomological degrees $s \leq 128$, 
\item 
\label{Abasis}
chain maps lifting each member in the resulting basis for
$\Ext_\cA^{s,t}(\bF_2,\bF_2)$ in this range, and 
\item
a chain map which gives the Hopf algebra squaring operation 
\[
Sq^0 : \Ext_\cA^{s,t}(\bF_2,\bF_2) \lra \Ext_\cA^{s,2t}(\bF_2,\bF_2) \,.
\]
\end{enumerate}

This document describes the files containing this information.  The
resolution was produced by the first author's software {\tt ext},
version 1.9.3.   This is contained in the file {\tt ext.1.9.3.tar.gz}.
The remaining contents of the top level directory are a copyright
notice,  a listing ({\tt ls-lR.txt}) of the files herein,
and a directory {\tt A}.
The latter has a subdirectory {\tt S-200} which contains
the data, and a subdirectory {\tt src} containing the source code for
a C program not included in {\tt ext.1.9.3.tar.gz}.  This program
writes a script to create cocycles for each basis element in~(\ref{Abasis})
above.

The resolution together with the chain maps give $\Ext_\cA(\bF_2,\bF_2)$
as an algebra.  It is worth pointing out that the minimal resolution
and the chain maps contain vastly more information than this.   Much
secondary and higher order structure is available from this data,
as well.   The discussion of Toda brackets (Massey products) below,
and the data contained in the file {\tt himults} are examples:
the file {\tt himults} contains information about products by 
elements of cohomological degree 1 obtained without use of chain maps, while
the Toda brackets files are produced from the chain maps without
reference to any null-homotopies.

The package {\tt ext} is designed to produce minimal resolutions
\[
0 \lla M \stackrel{d_0}{\lla} C_0 \stackrel{d_1}{\lla} 
        C_1 \stackrel{d_2}{\lla} \cdots \stackrel{d_{S}}{\lla} C_{S}
\]
and
\[
0 \lla N \stackrel{d_0}{\lla} D_0 \stackrel{d_1}{\lla} 
        D_1 \stackrel{d_2}{\lla} \cdots \stackrel{d_{S}}{\lla} D_{S}
\]
for any finite $\cA$-modules $M$ and $N$, and to lift cocycles $x
\in \Ext_\cA^{s_0,t_0}(M,N)$ represented as $\cA$-homomorphisms
$x : C_{s_0} \lra \Sigma^{t_0} N$ to chain maps $\{C_{s_0+s} \lra
\Sigma^{t_0} D_{s}\}_{s}$.  The exposition here is focused on the
case $M = N = \bF_2$, and the range $0 \leq s \leq 128$, $0 \leq t
\leq 200$, but at various points it may be useful to remember the
extra generality of the code.

\section{The resolution}
\label{sec:resolution}

Let us write
\[
0 \lla \bF_2 \stackrel{d_0}{\lla} C_0 \stackrel{d_1}{\lla} 
        C_1 \stackrel{d_2}{\lla} \cdots \stackrel{d_{128}}{\lla} C_{128}
\]
for our minimal resolution.  It is significant that minimality allows us to
identify $\Hom_\cA^t(C_s, \bF_2)  = \Hom_\cA(C_s, \Sigma^t \bF_2)$ with
$\Ext_\cA^{s,t}(\bF_2, \bF_2)$.

Since our resolution is limited to the range $0 \leq
s \leq 128$ and $0 \leq t \leq 200$, it is simply an initial segment
of a resolution.   For brevity, we shall nonetheless call it ``the
resolution'' and  write $\Ext$ for that part of
$\Ext_\cA(\bF_2,\bF_2)$ which lies in this range.
The resolution is described by the following files.

\begin{enumerate}
\item
{\tt Def}.  This file defines the $\cA$-module $\bF_2$ as the module
which is $1$-dimensional over $\bF_2$ with its sole generator in
degree $0$.

\item
{\tt MAXFILT}.  This contains the maximum cohomological degree, $S=128$,
through which the resolution is calculated.

\item
{\tt Shape}.  This file describes the free $\cA$-modules $C_s$.   Its
first entry, $128$, gives the maximum cohomological degree $S$.  This
is followed by the $129$ integers $\dim_{\cA}(C_s)$ for $s = 0, 1,
\ldots, 128$.  This is followed by the internal degree of each of
these generators, first for $C_0$, then for $C_1$, up to $C_{128}$.
This data determines $\Ext$ as a bigraded $\bF_2$-vector space.

We write {\tt s\_g} or $s_g$ for the cocycle dual to the {\em
$0$-indexed} $g^{\mathrm{th}}$ generator of $C_s$.   Thus, $0_0 \in
\Ext^{0,0} = \Ext_\cA^{0,0}(\bF_2,\bF_2)$ is the unit $1 : C_0 \lra \bF_2$ dual
to the $\cA$-module generator of $C_0$, while $1_i \in
\Ext^{1,2^i} = \Ext_\cA^{1,2^i}(\bF_2,\bF_2)$ is the ``Hopf map'' $h_i : C_1 \lra
\Sigma^{2^i}\bF_2$ dual to the $\cA$-module generator of $C_1$ that
$d_1$ sends to $Sq^{2^i}$.  When we need to refer to the $\cA$-module
generators of $C_s$, we shall write them as $s_g^*$ or {\tt s\_g*}.

\item
\label{item:elementformat}
{\tt Diff.s} and {\tt hDiff.s} for each $s$, $0 \leq s \leq 128$.
These files contain the differentials $d_s$ for $0 \leq s \leq 128$.
The file {\tt Diff.s} contains the elements $d_s(s_g^*)$, for
$g=0,1,\ldots$ in that order.  Each is preceded by its internal
degree.  The file {\tt Diff.s} is written in a condensed format
which takes the least space possible for the formats legible to the
program {\tt ext}.   The files {\tt hDiff.s} are provided as ``humanly
readable Diff files'', and write the differentials using the Milnor
basis.   A couple of examples should suffice to show how to read
them.  The file {\tt hDiff.1} starts

\begin{verbatim}
         8        200
1

1
0 1 1 i(1).

2

1
0 2 1 i(2).

4

1
0 4 2 i(4).
...
\end{verbatim}
which should be read as saying that the $\cA$-dimension of $C_1$
is $8$, that the part of the resolution given here is complete through
internal degree~$200$, and that
\begin{enumerate}
\item
generator $1_0^*$ lies in degree 1, and $d_1(1_0^*) = Sq^1(0_0^*)$,
\item
generator $1_1^*$ lies in degree 2, and $d_1(1_1^*) = Sq^2(0_0^*)$,
\item
generator $1_2^*$ lies in degree 4, and $d_1(1_2^*) = Sq^4(0_0^*)$,
et cetera.
\end{enumerate}
In the file {\tt hDiff.2}, the fifth entry, i.e., the entry for 
$d_2(2_4^*)$, is 
\begin{verbatim}
9

3
0 8 4 i(8)(2,2).
1 7 4 i(7)(4,1)(0,0,1).
3 1 1 i(1).
\end{verbatim}
This says that $2_4^*$ has internal degree $9$, and that
\begin{align*}
d_2(2_4^*) &= (Sq^8 + Sq^{(2,2)})(1_0^*)  \\
   &\quad + (Sq^7 + Sq^{(4,1)} + Sq^{(0,0,1)})(1_1^*)  \\
   &\quad +  Sq^1(1_3^*) \,.
\end{align*}
Here, the initial $3$  in the description of $d_2(2_4^*)$ says that
$d_2(2_4^*)$ is the sum of three terms, and the subsequent lines
describe those terms.  The first line, ``{\tt 0 8 4 i(8)(2,2).}'',
tells us that the first term is a multiple of $1_0^*$ with the
degree 8 coefficient $Sq^8 + Sq^{(2,2)}$.   The number $4$ here is
the $\bF_2$-dimension of $\cA$ in degree 8, and is used by the
program {\tt ext} to determine the amount of space which must be
allocated.  Note that elements of the Steenrod algebra are written
in the Milnor basis, not the admissible basis.\footnote{ 
The small ``i'' indicates that the notation is {\em
internal} to the programs defining the Steenrod algebra.  The main
body of the code will compute minimal resolutions for any connected
augmented $\bF_2$-algebra, and the other notations for coefficients
are generic notations for bitstrings which are independent of the
algebra (see Section~\ref{sec:chain} for some discussion of formats
``x'' and ``s'').}

\begin{figure}
\includegraphics{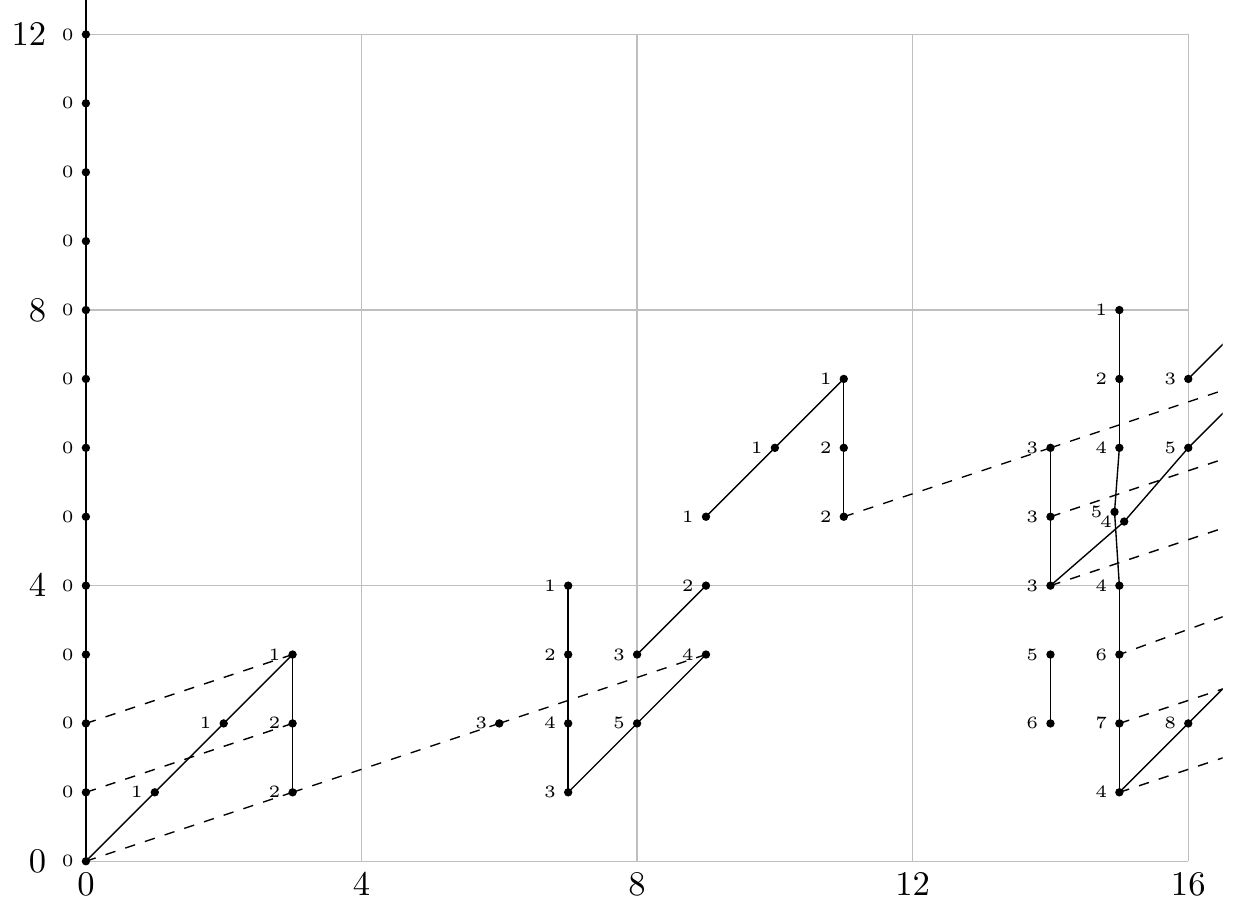}
\caption{$\Ext_\cA^{s,n+s}(\bF_2, \bF_2)$,\,\,
$0 \leq n \leq 16,\,\, 0 \leq s \leq 12$}
\label{fig:samplechart1}
\end{figure}

\begin{figure}
\includegraphics{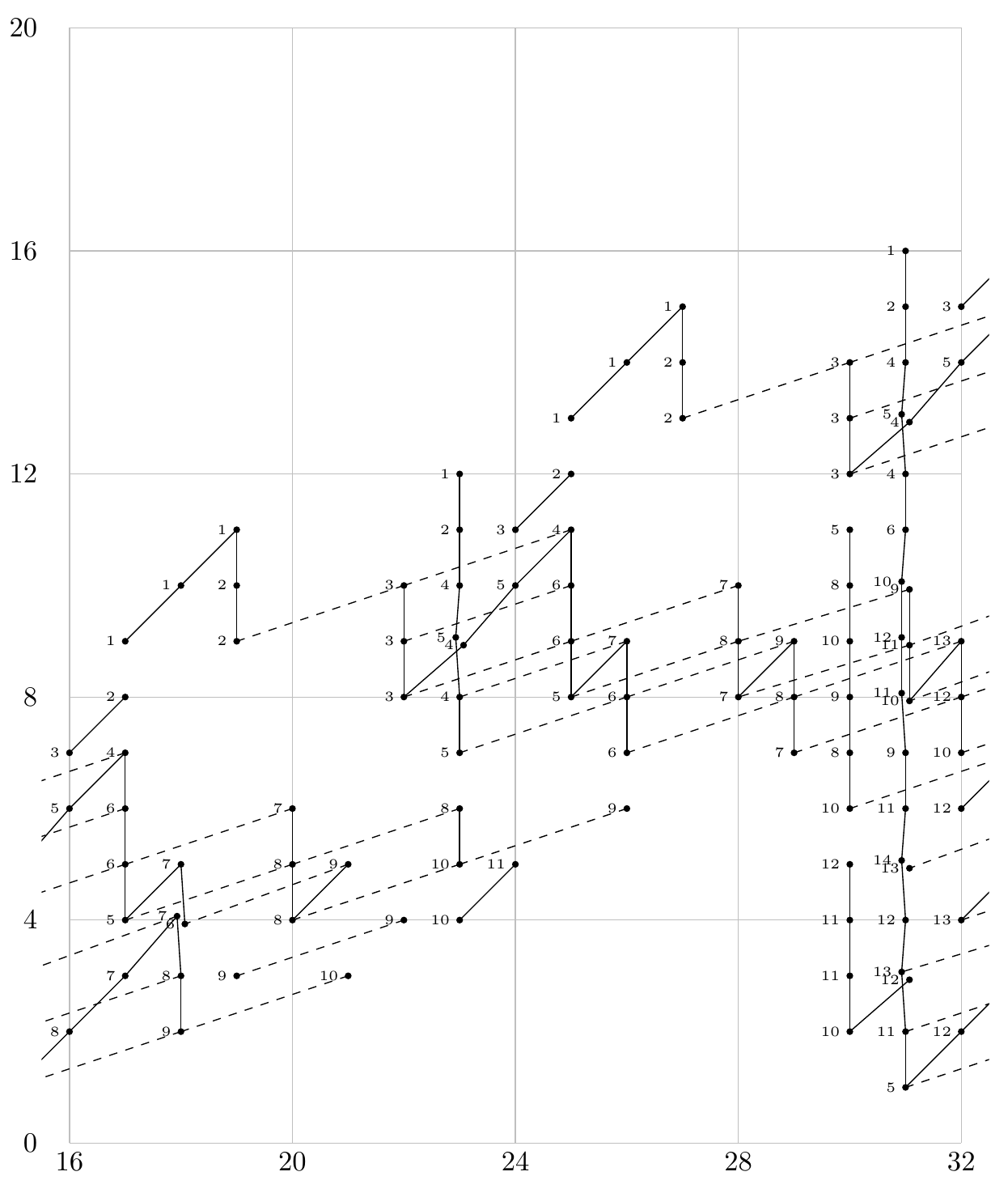}
\caption{$\Ext_\cA^{s,n+s}(\bF_2, \bF_2)$,\,\,
$16 \leq n \leq 32,\,\, 0 \leq s \leq 20$}
\label{fig:samplechart2}
\end{figure}

\begin{figure}
\includegraphics{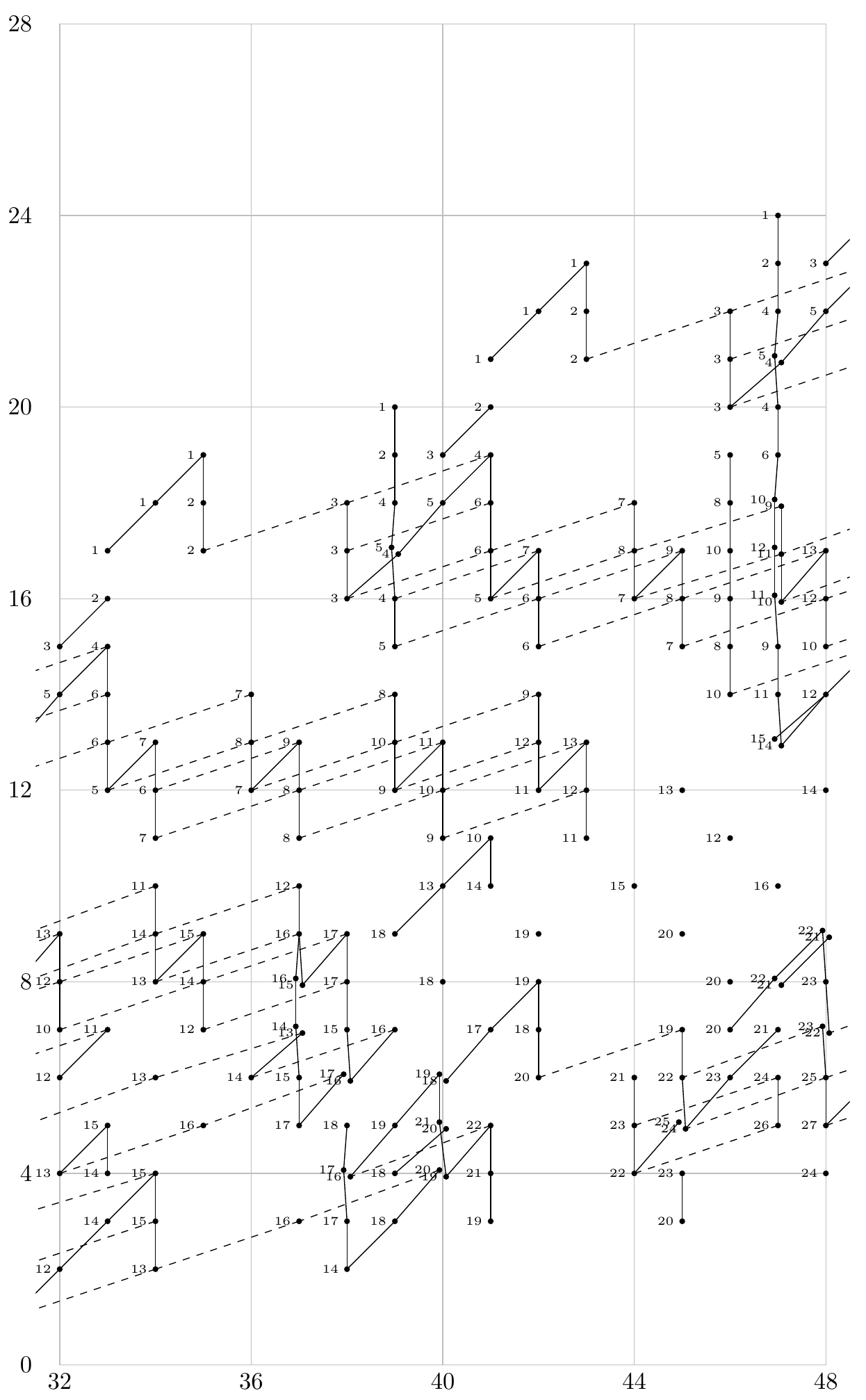}
\caption{$\Ext_\cA^{s,n+s}(\bF_2, \bF_2)$,\,\,
$32 \leq n \leq 48,\,\, 0 \leq s \leq 28$}
\label{fig:samplechart3}
\end{figure}


\item
{\tt Ext-A-F2-F2-0-200.tex} and {\tt Ext-A-F2-F2-0-200.pdf}.
These figures are in the usual ``Adams chart'' format
showing the resolution together with the action of the
elements $h_0$, $h_1$ and $h_2$.  
The program {\tt chart} included in {\tt ext.1.9.3} allows
one to make visualizations like this for any box
$s_{lo} \leq s \leq s_{hi}$, $t_{lo} \leq t \leq t_{hi}$.
It is also possible to simply use the {\tt TikZ} command
{\tt clip} to extract desired sections from the \TeX{} file
we provide.

In these charts, each cocycle $s_g$ is represented by a filled
circle with the sequence number, $g$, written to its left, as in
Figures~\ref{fig:samplechart1} to~\ref{fig:samplechart3}.  For
example, in bidegree $(n,s) = (t-s,s) = (15,5)$ we see that $5_4 = h_1 \cdot
4_3$ and $5_5 = h_0 \cdot 4_4$.
See also Figure~\ref{Ext-A-F2-F2-0-200+gray} at the end of this document,
which shows the full Adams chart in the range $t \le 200$.

\item
{\tt S-200.tex} and {\tt S-200.pdf}.
These give a stem-by-stem list of the results together with
the products by the $h_i$.

\item
{\tt Maxt} and {\tt himults}.  These are generated by the
program {\tt report} included in {\tt ext.1.9.3} and are
used by the {\tt chart} and {\tt vsumm} commands to create
the Adams charts and the \TeX~ summary.
The file {\tt Maxt} records, for each cohomological degree $s$,
the degree $t$ through which {\tt Diff.s} is complete.  The file
{\tt himults} records the $h_i$ products using the observation
cited in the first paragraph of Section \ref{sec:brackets} to
extract this from the {\tt Diff.s} files.   The format of
each entry is
\begin{verbatim}
s g s0 g0 i
\end{verbatim}
meaning that the product $h_i \cdot (s_0)_{g_0}$ contains the term $s_g$.
For example, in bidegree $(n,s) = (t-s,s) = (78,8)$ we have three
generators $8_{60}$, $8_{61}$ and~$8_{62}$, and the thirteen lines
\begin{verbatim}
8 60 7 33 4
8 60 7 34 4
8 60 7 56 1
8 60 7 58 1
8 61 7 33 4
8 61 7 34 4
8 61 7 48 3
8 61 7 56 1
8 61 7 57 1
8 61 7 60 0
8 62 7 2 6
8 62 7 35 4
8 62 7 61 0
\end{verbatim}
say that
\begin{align*}
h_0 \cdot 7_{60} & = 8_{61} \\
h_0 \cdot 7_{61} & = 8_{62} \\
h_1 \cdot 7_{56} & = 8_{60} + 8_{61} \\
h_1 \cdot 7_{57} & = 8_{61} \\
h_1 \cdot 7_{58} & = 8_{60} \\
h_3 \cdot 7_{48} & = 8_{61} \\
h_4 \cdot 7_{33} & = 8_{60} + 8_{61} \\
h_4 \cdot 7_{34} & = 8_{60} + 8_{61} \\
h_4 \cdot 7_{35} & = 8_{62} \\
h_6 \cdot 7_{2} & = 8_{62} \,,
\end{align*}
while the other products $h_i \cdot (s_0)_{g_0}$ in this bidegree
are zero.  This is used by {\tt chart} to draw the lines representing
$h_0$, $h_1$ and $h_2$ products.  
\end{enumerate}

\section{Products}
\label{sec:products}

We compute products in $\Ext$ by composing chain maps.  Suppose
the classes $x \in \Ext_\cA^{s_0,t_0}(N,P)$ and
$y \in \Ext_\cA^{s_1,t_1}(M,N)$ are represented by cocycles
\[
x : D_{s_0} \lra \Sigma^{t_0} P \qquad {\mathrm{ and}} \qquad
y : C_{s_1} \lra \Sigma^{t_1} N.
\]
Then $\Sigma^{t_1} x \circ y_{s_0}$ is a cocycle representing the
product $xy$, where $\{y_s\}_s$ is a chain map lifting $y$.  The
chain map with components $\Sigma^{t_1} x_s \circ y_{s+s_0}$ is a
lift of this cocycle.

\[
\xymatrix@C=3.1ex{
M
&
C_0
\ar[l]
&
{}\phantom{x}
\cdots
\phantom{x}
\ar[l]
&
C_{s_1}
\ar_(0.6){y}[dl]
\ar^{y_0}[d]
\ar[l]
&
{}\phantom{x}
\cdots 
\phantom{x}
\ar[l]
&
C_{s_0+s_1}
\ar^{y_{s_0}}[d]
\ar[l]
&
{}\phantom{x}
\cdots
\phantom{x}
\ar[l]
&
C_{s+s_0+s_1}
\ar_{y_{s+s_0}}[d]
\ar[l]
\\
&&
\Sigma^{t_1} N
&
\Sigma^{t_1} D_0
\ar[l]
&
{}\phantom{x}
\cdots
\phantom{x}
\ar[l]
&
\Sigma^{t_1} D_{s_0}
\ar_(0.6){\Sigma^{t_1}x}[dl]
\ar^{\Sigma^{t_1}x_0}[d]
\ar[l]
&
{}\phantom{x}
\cdots 
\phantom{x}
\ar[l]
&
\Sigma^{t_1} D_{s+s_0}
\ar_{\Sigma^{t_1}x_s}[d]
\ar[l]
\\
&&&&
\Sigma^{t_0+t_1} P
&
\Sigma^{t_0+t_1} E_0
\ar[l]
&
{}\phantom{x}
\cdots
\phantom{x}
\ar[l]
&
\Sigma^{t_0+t_1} E_{s}
\ar[l]
\\
}
\]

In our situation, $M = N = P = \bF_2$ and $E_s = D_s = C_s$.  In
this case, the composite $\Sigma^{t_1} x \circ y_{s_0}$ is determined
by recording those generators $(s_0+s_1)_{g}^*$ of $C_{s_0+s_1}$
that are mapped nontrivially.    When the cocycle $x$ is $(s_0)_{g_0}$,
dual to a generator $(s_0)_{g_0}^*$ of $C_{s_0}=D_{s_0}$, the list
of such generators is the set of those whose image under $y_{s_0}$
contains a term $1 \cdot (s_0)_{g_0}^*$, since all other terms will
be sent to $0 \in \bF_2$ by $(s_0)_{g_0}$.   The program {\tt
collect} in the {\tt ext} package gleans this information from the
chain map files described in the next section and organizes it into
the file {\tt A/S-200/all.products}.

Each line in the file {\tt all.products} has the form
\begin{verbatim}
         s   g  (  s0    g0     F2)  s1_g1
\end{verbatim}
which means that the chain map lifting the cocycle $(s_1)_{g_1}$, applied to the basis
element $s_g^*$, contains the term $1 \cdot (s_0)_{g_0}^*$.   

\begin{proposition}
$(s_0)_{g_0} \cdot (s_1)_{g_1}$ is the sum of all such $s_g$.
\end{proposition}

The file {\tt all.products} is organized so that each paragraph
lists, for a cocycle $s_g$, all such pairs $(s_0)_{g_0}$, $(s_1)_{g_1}$
with the $(s_0)_{g_0}$ in increasing order.\footnote{To be precise,
the new program {\tt collect.JR}, included in the directory {\tt
S-200}, sorts the entries in this manner.   It will replace the old
{\tt collect} in the next release of the {\tt ext} package.
References to {\tt collect} here should generally be interpreted to
refer to the improved {\tt collect.JR}.}

To compute $(s_0)_{g_0} \cdot (s_1)_{g_1}$ it is best to start by
noting which $s_g$ span the bidegree containing the product, so as
not to miss a term.  For example, consider the products landing in
$\Ext^{7,7+37}$.  This bidegree is $2$-dimensional, spanned by
$7_{13}$ and $7_{14}$.  The entries in the file {\tt all.products} for these
two cocycles are 
\begin{verbatim}
  7   13  (  0    0     F2)  7_13
  7   13  (  1    1     F2)  6_14
  7   13  (  1    2     F2)  6_13
  7   13  (  1    3     F2)  6_10
  7   13  (  2    3     F2)  5_13
  7   13  (  3    9     F2)  4_6
  7   13  (  4    6     F2)  3_9
  7   13  (  5   13     F2)  2_3
  7   13  (  6   10     F2)  1_3
  7   13  (  6   13     F2)  1_2
  7   13  (  6   14     F2)  1_1

  7   14  (  0    0     F2)  7_14
  7   14  (  1    0     F2)  6_15
  7   14  (  1    3     F2)  6_10
  7   14  (  2    0     F2)  5_17
  7   14  (  5   17     F2)  2_0
  7   14  (  6   10     F2)  1_3
  7   14  (  6   15     F2)  1_0
\end{verbatim}
This says that
\begin{align*}
7_{13} & = 0_0 \cdot 7_{13} 
       = 1_1 \cdot 6_{14}
       = 1_2 \cdot 6_{13}
       = 2_3 \cdot 5_{13}
       = 3_9 \cdot 4_{6}\\
7_{14} & = 0_0 \cdot 7_{14}
       = 1_0 \cdot 6_{15}
       = 2_0 \cdot 5_{17}\\
7_{13} + 7_{14} & = 1_3 \cdot 6_{10} \,.
\end{align*}
Eliminating redundant cases and using the traditional notations
$h_i = 1_i$, $h_0^s = s_0$, $t=6_{14}$, $n=5_{13}$,  $r=6_{10}$,
$c_1 = 3_9$, $f_0 = 4_6$, $x=5_{17}$, these two paragraphs say that
the  two elements $h_1 t = h_2^2 n = c_1 f_0$ and $h_0^2 x$ span
$\Ext^{7,7+37}$, that $h_3 r = h_1 t + h_0^2 x$, and that all other
products landing in this bidegree are $0$.  The $h_0, h_1$ and $h_2$
products in this bidegree can be seen in the chart shown in
Figure~\ref{fig:samplechart3}.   (The original definitions of $f_0$
did not distinguish between $f_0$ and $f_0 + h_1^3 h_4$. We eliminate
this ambiguity by defining $f_0 = Sq^1(c_0)$.  It is shown in
\cite{BNT} that $Sq^1(c_0) = 4_6$.)

\begin{remark}
\label{moreindec}
Cocycles $s_g$ whose paragraph in {\tt all.products} has only a single entry,
$0_0 \cdot s_g$, are clearly indecomposable.  There are~$1167$ of these
in the range calculated, lying in filtrations~$1$ through~$65$.  However,
to identify all the indecomposables, simple textual work is
insufficient:  we must calculate $\fm/\fm^2$, where $\fm$ is the
maximal ideal of $\Ext$.  For example $\Ext^{7,7+133}$ is spanned
by $7_{124}$, $7_{125}$ and $7_{126}$ with products $h_3 \cdot
6_{97} = 7_{124}+7_{125}$ and $h_1 \cdot 6_{107} = 7_{126}$.
Thus, $7_{124}$ and $7_{125}$ project to the same nonzero element
of $\fm/\fm^2$.   Similarly, $h_4 \cdot 8_{114} = 9_{178} + 9_{179}$
with each term indecomposable, but equal to one another modulo
decomposables.
\end{remark}

\section{Chain maps}
\label{sec:chain}

A cocycle $s_g$ of internal degree $t$ is an $\cA$-module homomorphism
$s_g : C_s \lra \Sigma^t \bF_2$.   Using {\tt ext} this can be
lifted to a chain map of bidegree $(s,t)$,
\[
\{C_{s+s_0} \lra \Sigma^t C_{s_0} \mid 0 \leq s_0 \leq s+s_0 \leq 128\}.
\]
The data describing the cocycle $s_g$ and
our chain map lifting the cocycle are in the subdirectory {\tt s\_g} of the
directory {\tt S-200} containing the resolution.    The relevant
files are as follows.

\begin{enumerate}
\item {\tt maps} and subdirectories {\tt A/S-200/s\_g}.
In {\tt A/S-200}, the file {\tt maps} is a list of all the cocycles $s_g$.
Each of these has a subdirectory {\tt A/S-200/s\_g} which contains the
cocycle's definition, the chain map lifting it, and data derived
from this.

In order to study various aspects of the resolution it may be useful
to focus on a smaller set of maps.   The package {\tt ext} has
utilities, such as {\tt collect}, which operate on such lists.  For
example, the invocation
\begin{verbatim}
./collect somemaps someproducts
\end{verbatim}
would create a file named {\tt someproducts} containing all products
by maps listed in {\tt somemaps}.

\item 
{\tt s\_g/Def}.
This file contains the definition of the cocycle $s_g$.  If
the internal degree of $s_g$ is $t$, the file {\tt s\_g/Def}
will contain
\begin{verbatim}
 s t F2 F2 s_g 1

g

1
0 0 1 x80
\end{verbatim}
The reader who simply wants to use this data will not need the
following discussion of cocycle definition files, but for 
completeness, here is the meaning of these entries.

Reading them in order, it says that $s_g : C_s \lra \Sigma^t \bF_2$
is stored in the subdirectory {\tt s\_g} and maps exactly one
generator of $C_s$ nontrivially, namely generator number $g$, sending
it to the unique nontrivial element of $\bF_2$, which is the $\bF_2$
basis element numbered $0$.   The format of a general {\em cochain}
$x : C_s \lra \Sigma^t N$, where $C_*$ is a resolution of $M$, is

\begin{verbatim}
 s t M N x n

g1
x(g1)

g2
x(g2)

...

gn
x(gn)
\end{verbatim}
This specifies the values in $N$ of $x$ on the $n$ generators
numbered $g_1$ through $g_n$.  Generators of $C_s$ which are not
mentioned are mapped to $0$.  This defines a {\em cocycle} iff it
can be lifted to a chain map iff it can be lifted over the first
stage, $d_1 : D_1 \lra D_0$, of a resolution of $N$.   Thus, the
process of lifting can be used to check that a cochain $x$ is a
cocycle.   By minimality of the resolutions we produce, all cochains
$x : C_s \lra \Sigma^t \bF_2$ are cocycles.

For example, if $N$ is $4$-dimensional over $\bF_2$, with basis
elements $0$ to $3$ in degrees $0$, $1$, $1$ and~$2$, respectively,
and a cochain $x : C_5 \lra \Sigma^{12} N$ had nonzero values $x(5_3)
= 1 + 2$ and $x(5_6) = 3$, the cochain definition file {\tt x/Def}
could be

\begin{verbatim}
 5 12 M N x 2

3

2
1 0 1 x80
2 0 1 s0.

6

1
3 0 1 i(0).
\end{verbatim}
Here, ``{\tt 0 1 x80}'', ``{\tt 0 1 s0.}'' and ``{\tt 0 1 i(0).}''
are three different ways to write the unit $1 \in \cA_0$.  The
initial ``{\tt 0 1}'' in each says that they describe an element
of $\cA_0$, which is $1$-dimensional over $\bF_2$.   The first,
``{\tt x80}'', uses hexadecimal notation to express the binary
vector with a $1$ in its first (and only) entry.   The second,
``{\tt s0.}'', lists the sequence number, $0$, of the coordinates
whose entries are $1$ rather than $0$.   The third, ``{\tt i(0).}'',
writes the operation in Milnor basis form: $Sq^0 = 1 \in \cA_0$.
We could (and usually do) write the same information in the form

\begin{verbatim}
 5 12 M N x 2

3

2
1 0 1 x80
2 0 1 x80

6

1
3 0 1 x80
\end{verbatim}
for simplicity and uniformity.

\item 
{\tt s\_g/Map} and {\tt s\_g/Map.aug}.  The {\tt ext} package writes
the chain map lifting the cocycle $s_g$ in the file {\tt s\_g/Map}.
This file is simply a list of entries of the form
\begin{verbatim}
s0 g0 0
\end{verbatim}
meaning that the lift of $s_g$ sends {\tt s0\_g0*} to $0$, or
\begin{verbatim}
s0 g0 
y
\end{verbatim}
where {\tt y} is the representation, in the manner discussed in
item~(\ref{item:elementformat}) of Section~\ref{sec:resolution}, of
the nonzero image of {\tt s0\_g0*} under the chain map lifting
$s_g$.   The entries in the {\tt Map} file do not have to occur in
any particular order.   There is a program, {\tt checkmap s}, which 
determines the elements, if any,  mapping to filtration $\leq s$
which are not yet present in the {\tt Map} file.   This has been used
to verify the completeness of the data we present here.

The file {\tt s\_g/Map.aug} is extracted from the {\tt Map} file
by applying the augmentation $\cA \lra \bF_2$, i.e., by discarding
all terms $a \cdot s_g^*$ in which $a\in\cA$ has positive degree.
The information in this file is used by {\tt collect} to compile
the file {\tt all.products} described in the previous section.

\item
{\tt s\_g/brackets} and 
{\tt s\_g/brackets.sym}.
These are discussed in the next section.
\end{enumerate}

\section{Toda brackets}
\label{sec:brackets}

Products $h_i \cdot s_g$ can be calculated directly from the
resolution without computing either of the chain maps lifting $h_i$
or $s_g$.  Precisely, $h_i \cdot s_g$ is the sum of those $(s{+}1)_{g_1}$
such that $d_{s+1}((s{+}1)_{g_1}^*)$ contains the term $Sq^{2^i}\cdot s_g^*$.
In a similar way, having computed the chain maps, we are now able
to evaluate all Toda brackets\footnote{We take the point of view
that secondary products with respect to composition should be called
{\em Toda brackets}, while secondary products in a DGA should be
called {\em Massey products}.  By the usual device of considering
the DGA of endomorphisms of a chain complex, the two are equivalent.}
of the form $\langle h_i, (s_0)_{g_0}, (s_1)_{g_1} \rangle$.

\begin{proposition}
\label{prop:brackets}
If $h_i \cdot (s_0)_{g_0} = 0$ and $(s_0)_{g_0} \cdot (s_1)_{g_1}
= 0$ then the Toda bracket $\langle h_i, (s_0)_{g_0}, (s_1)_{g_1}
\rangle$ contains the sum of all those $s_g$ such that the chain
map lifting $(s_1)_{g_1}$ applied to $s_g^*$ contains a term
$Sq^{2^i} \cdot (s_0)_{g_0}^*$.
\end{proposition}

Conceptually, products by elements of cohomological degree 1 are visible in the
resolution reduced modulo $\fm^2$, where $\fm \subset \cA$ is the
augmentation ideal, while brackets with first entry of cohomological
degree 1 are visible in the chain maps reduced modulo $\fm^2$.   This data
is extracted from the {\tt Map} file and placed in the files
{\tt brackets} and {\tt brackets.sym}.  The files contain the same information,
but {\tt brackets.sym} is easier for a human to read.   Each entry in
{\tt s1\_g1/brackets.sym} will have the form
\begin{verbatim}
s_g in < hi, g0, s1_g1 >
\end{verbatim}
Since the filtration $s$ of the bracket is $1 + s_0 + s_1 -1 = s_0 + s_1$,
we deduce that the middle entry is {\tt s0\_g0} with
$s_0 = s-s_1$.

As with products, it is important to remember that the value of the bracket
is the sum of all the $s_g$ which appear, so that it is necessary to survey
all the possible terms in the relevant bidegree before reaching conclusions.

We also hold the point of view that it is sufficient to produce one
element of the bracket, with the other elements being obvious from
the known indeterminacy.  We also note that the {\tt brackets.sym}
file is simply recording information about the chain map, and that
it is the responsibility of the user to check whether the bracket
is defined.

For example, the first few entries in {\tt 1\_0/brackets.sym}  are
\begin{verbatim}
2_23 in < h7, 0, 1_0 >
2_8 in < h4, 0, 1_0 >
2_1 in < h0, 1, 1_0 >
2_5 in < h3, 0, 1_0 >
3_3 in < h3, 0, 1_0 >
3_3 in < h1, 3, 1_0 >
\end{verbatim}
Since $h_0 h_7 \neq 0$ and $h_0 h_4 \neq 0$, the first two brackets
are not defined.  The third says that $h_1^2 = 2_1 \in \langle h_0,
h_1, h_0\rangle$, a familiar consequence of Hirsch's formula $y(x
\cup_1 x) \in \langle x,y,x \rangle$.  The next two are also not
defined.   Following that we find $c_0 = 3_3 \in \langle h_1, h_2^2,
h_0\rangle$.

\section{$Sq^0$}

In general, the cocommutative Hopf algebra Steenrod operations 
\[
Sq^i : \Ext_\cA^{s,t}(\bF_2,\bF_2) \lra \Ext_\cA^{s+i,2t}(\bF_2,\bF_2)
\]
are computationally quite expensive
(\cite{BruYoneda}*{Problem Session} and \cite{BNT}).   However, the two extremes
$Sq^s : \Ext^{s,t} \to \Ext^{2s,2t}$ and $Sq^0 : \Ext^{s,t} \to \Ext^{s,2t}$ are 
easily calculated using {\tt ext}.   The first is simply the squaring
operation $Sq^s(x) = x^2$ for $x \in \Ext^{s,t}$, which we have already discussed.
At the other extreme,
in \cite{May70}*{Proposition 11.10}, it is
shown that the operation $Sq^0$ can be calculated by $Sq^0([a_1 | \ldots
| a_s]) = [a_1^2 | \ldots | a_s^2]$ in the cobar complex
for the dual Steenrod algebra.  This implies that if $\Phi \cA_*$
is the double of the dual Steenrod algebra, in which the degrees of all
the elements are doubled, then $Sq^0 \: \Ext^{s,t}_{\cA_*}(\bF_2, \bF_2) \to
\Ext^{s,2t}_{\cA_*}(\bF_2, \bF_2)$ is induced by the degree-preserving Hopf
algebra homomorphism $F \: \Phi \cA_* \to \cA_*$
that sends $\Phi\xi_i$ to $\xi_i^2$ for each $i\ge1$.  Dually, it
is induced by the degree-preserving Hopf algebra homomorphism $V \: \cA
\lra \Phi \cA$ that sends an ``even''
Milnor basis element $Sq^{(2r_1, \ldots, 2r_k)}$
to $\Phi Sq^{(r_1,\ldots, r_k)}$, and other Milnor basis elements to $0$.
Restricting along this homomorphism gives
\[
Sq^0 \: \Ext^{s,t}_\cA(\bF_2,\bF_2) \cong \Ext^{s,2t}_{\Phi \cA}(\bF_2,\bF_2)
        \lra \Ext^{s,2t}_\cA(\bF_2,\bF_2).
\]
A slight modification of the computer code that calculates chain maps can compute this:  a program
{\tt startsq0} computes the restriction
$V_{s-1}(d(x))$ for each generator $x = s_g^*$ in the minimal
$\cA$-module resolution $(C_*, d)$ of~$\bF_2$, and the same program
that computes lifts for chain maps then solves for an element $V_s(x)$
satisfying $d(V_s(x)) = V_{s-1}(d(x))$.  We recover $Sq^0$
as $\Hom_\cA(V_*, \bF_2)$.  This inductive calculation is begun by setting
$V_0(0_0^*) = 0_0^*$, so that $Sq^0(1) = 1$.  (This discussion is quoted from
our proof of Proposition 11.26 in~\cite{BR21}.)

The files involved are
\begin{enumerate}
\item
{\tt dosq0} and {\tt maps.sq0}.  
The first is a script to run the computation of the lifts, while the
second tells which subdirectory contains the map $Sq^0$.
\item 
{\tt S-200/Sq0/Map} and {\tt S-200/Sq0/Map.aug}.
This contains the data defining the chain map $V$  and its reduction
modulo the augmentation ideal $\fm \subset \cA$, respectively. This
latter defines the dual of the map $Sq^0$.
\item
{\tt all.sq0}.
This contains the data in the {\tt Sq0/Map.aug} file in the format of the {\tt all.products} file.
For example, its entries
\begin{verbatim}
  0    0  (  0    0     F2)  Sq0
  1    1  (  1    0     F2)  Sq0
  1    2  (  1    1     F2)  Sq0
  ...
  3    9  (  3    3     F2)  Sq0
  3   19  (  3    9     F2)  Sq0
  3   34  (  3   19     F2)  Sq0
  3   55  (  3   34     F2)  Sq0
  ...
  4   16  (  4    5     F2)  Sq0
  4   19  (  4    6     F2)  Sq0
\end{verbatim}
show that 
\begin{align*}
Sq^0(1)  = Sq^0(0_0) & = 0_0  = 1 \\
Sq^0(h_0)  = Sq^0(1_0) & = 1_1   = h_1 \\
Sq^0(h_1)  = Sq^0(1_1) & = 1_2  = h_2 \\
\ldots \\
Sq^0(c_0)  = Sq^0(3_3) & = 3_9  = c_1 \\
Sq^0(c_1)  = Sq^0(3_9) & = 3_{19}  = c_2  \\
Sq^0(c_2)  = Sq^0(3_{19}) & = 3_{34}  = c_3 \\
Sq^0(c_3)  = Sq^0(3_{34}) & = 3_{55}  = c_4 \\
\ldots \\
Sq^0(e_0)  = Sq^0(4_5) & = 4_{16}  = e_1 \\
Sq^0(f_0)  = Sq^0(4_6) & = 4_{19}  = f_1 \,.
\end{align*}
Later we find
\begin{verbatim}
  6  102  (  6   33     F2)  Sq0
  6  103  (  6   32     F2)  Sq0
  6  103  (  6   33     F2)  Sq0
  6  104  (  6   32     F2)  Sq0
\end{verbatim}
which reminds us that these files are the dual, i.e., chain level, data.
In $\Ext$, these say
\begin{align*}
Sq^0(A_0+A_0')  = Sq^0(6_{32}) & = 6_{103} + 6_{104}, \\
Sq^0(A_0)  = Sq^0(6_{33}) & = 6_{102} + 6_{103}, \,\,\,{\mathrm{and}}\,\,{\mathrm{hence}}\\
Sq^0(A_0')  = Sq^0(6_{32} + 6_{33}) & = 6_{102} + 6_{104} \,.
\end{align*}
\end{enumerate}

This information is used in Section~\ref{sec:concordance} to organize $\Ext$
into ``families'' linked by~$Sq^0$.

\section[A canonical basis]{A canonical basis for $\Ext_\cA(\bF_2,\bF_2)$}

The traditional, and by now familiar, notation for the elements
of $\Ext$ starts in a systematic fashion, with the elements
\[
h_i, P^k h_1, P^k h_2, c_i, P^k c_0, d_i, P^k d_0, e_i, P^k e_0, f_i,
g_i, i, j, k, \ldots \,,
\]
but becomes somewhat chaotic as the calculations are extended into 
higher bidegrees.   In the next section, we propose some ways of extending
this notation in a methodical fashion, but they do not suffice to give
names to elements of a basis for $\Ext$ even in the range we consider here,
let alone for all of $\Ext$.

In contrast, our $s_g$ form a well-defined canonical basis for $\Ext$,
which we now describe.  First, we totally order the terms $Sq^{R} s_g^*$ 
of $C_{s,t}$
by
\[
Sq^{R} s_g^* < Sq^{R'} s_{g'}^*
\]
iff
\begin{enumerate}
\item $g<g'$, or
\item $g=g'$ and $Sq^{R} < Sq^{R'}$, where the Milnor
basis elements $Sq^{R}$ are given reverse lexicographic order:
$(r_1, r_2, \ldots) < (r_1', r_2', \ldots)$ iff
for some $k$,
$r_k < r_k'$ and $r_i = r_i'$ for all $i>k$.
Thus,
\begin{multline*}
(n) < (n-3,1) < (n-6,2) < \cdots  \\
 < (n-7,0,1) < (n-10,1,1) < \cdots < (n-14,0,2) < (n-17,1,2) < \cdots  \\
 < (n-15,0,0,1) < (n-18,1,0,1) < \cdots < (n-22,0,1,1) < \cdots
\end{multline*}
\end{enumerate}

Each nonzero element $x \in C_{s,t}$ then has a {\em leading term}
$LT(x)$, which is the lowest term in $x$.

In the totally ordered basis $\{Sq^{R} s_g^*\}$ of a given bidegree
$(s,t)$, the decomposable elements, those with $\deg(Sq^{R}) > 0$,
form an initial segment which is followed by the
generators $s_g^*$ of bidegree $(s,t)$.

%

We can now inductively define our canonical basis as follows. 
We start with the bases $\{0_0^*\}$ for $C_0$ and $\{\}$ for
$C_s$ with $s>0$.  We may inductively assume given the basis
for $C_s$ in degrees less than $t$, and for $C_{s-1}$ in degrees
less than or equal to $t$.

{\bf Step 1:} Generating the image and kernel.\\
$\Im_{s,t}$ will be a totally ordered list of pairs $(x,dx)$ with
the leading terms of the $dx$ in strictly increasing order.
$\Ker_{s,t}$ will be a list of terms $x$.   Both are initially empty.

Consider the terms $Sq^{R} s_g^*$ in order.  Let
$x = Sq^{R} s_g^*$ and compute $dx = Sq^{R} d(s_g^*)$.   
Then, while $dx \neq 0$, if $LT(dx) = LT(dy)$ for a pair 
$(y,dy) \in \Im_{s,t}$, replace $x$ by $x-y$ and $dx$ by 
$dx -dy$.  If not, add $(x,dx)$ to $\Im_{s,t}$ and proceed to
the next decomposable term.  If instead $dx = 0$, add $x$
to the end of the list $\Ker_{s,t}$.

Note that the leading term of $dx$ will be increased
each time we replace $dx$ by $dx - dy$ until $dx$ either becomes
$0$ or has a leading term not already found among the $dy$ in
$\Im_{s,t}$.


{\bf Step 2:} Adding new generators.\\
We may inductively assume given $\Ker_{s-1,t}$.  For each
$x \in \Ker_{s-1,t}$, in order, let $c=x$.   Then, while $LT(x) = LT(dy)$
for some pair $(y,dy) \in \Im_{s,t}$, replace $x$ by $x-dy$.
If this process terminates with $LT(x) \neq 0$, add a new
generator $s_g^*$ with $d(s_g^*) = c$, then add a new pair 
$(z,x)$ to $\Im_{s,t}$, where $z$ is the difference of $s_g^*$ and those $y$
whose images $dy$ were subtracted from $c$ to get the final $x$
with a new leading term.
If the process terminates with $x=0$, do nothing.

\begin{remark}
We could choose, at this second step, to let $d(s_g^*) = x$
and add the pair $(s_g^*,x)$ to $\Im_{s,t}$.  The {\tt ext} code
prior to the year 2000 used that algorithm.
Experience shows
that the bases $s_g$ obtained from the algorithm described here
have $h_i \cdot s_g$ monomial far more frequently than those
produced by the old algorithm.
The Wayne State Research Report~\cite{Bru97} used the older
algorithm.  The first difference visible in $\Ext$ charts lies
in bidegree $(9,9+23)$.   In the new algorithm, $h_1 \cdot 8_3 = 9_4$
and $h_0 \cdot 8_4 = 9_5$.  In the older algorithm,
$h_0 \cdot 8_4 = 9_5$ also, but
$h_1 \cdot 8_3 = 9_4 + 9_5$.

The change alters the resolution much earlier, as can be seen by
doing hand calculations in low degrees.  In the old algorithm,
$d(2_1^*) = Sq^{(0,1)} \cdot 1_0^* + Sq^2 \cdot 1_1^*$, while the
new algorithm gives $d(2_1^*) = Sq^3 \cdot 1_0^* + Sq^2 \cdot 1_1^*$.
\end{remark}

\section{Concordance}
\label{sec:concordance}

In this section, we present the relation between our $s_g$ basis
and the notation used by other works on the cohomology of the
Steenrod algebra.   In the process, we make a natural extension
to the traditional notation using $Sq^0$.

The existing names are based on Tangora's calculation of the
$E_\infty$ term of the May spectral sequence in~\cite{Tan70} and
on Chen's Lambda algebra computation of $\Ext^s = \Ext_{\cA}^s(\bF_2, \bF_2)$ for $s \leq 5$
in~\cite{Che11}.   

There is some indeterminacy in the translation between our names
and both the May spectral sequence names and Chen's Lambda algebra
names.  

In the case of the May spectral sequence, elements of the $E_\infty$-term
of the May spectral sequence only determine elements of $\Ext$ up
to classes of higher May filtration.   In Table~\ref{concordance} we
indicate this by giving the indeterminacy of the May spectral
sequence name in parentheses.  For example, in bidegree $(5,5+62)$,
the May spectral sequence definition of $H_1(0)$ (written $H_1$ in
\cite{Tan70}) defines the coset $5_{32} + \langle 5_{33}, 5_{34}
\rangle$.  We denote this, in our table, by writing $32(33,34)$ in
the column giving the sequence number of the element.

The indeterminacy in relating Chen's Lambda algebra element names
to ours stems from the lack of a direct comparison between the two
complexes.  There are five such cases which we discuss in
Remark~\ref{Chentomin}.  Except for $T_0$, which is beyond the range
of Tangora's May spectral sequence calculation, this indeterminacy
is the same as that due to the May filtration.   In the case
of $T_0$, the indeterminacy is entirely due to the lack of a direct
comparison between the Lambda algebra and our minimal resolution.
For $T_0$, $Q_3(0)$ and $H_1(0)$, the indeterminacy could be
eliminated if we knew certain $h_i$-multiples, as noted in
Remark~\ref{Chentomin}.

The preprint \cite{IWX} of Isaksen-Wang-Xu studies the Adams
spectral sequence through $t-s \leq 95$.   They adopt some of
Tangora's notation for $\Ext$, augmented and regularized by the use
of operators which they call $M$, $\Delta$ and $\Delta_1$.   Their operator
$Mx$ is principally defined to be $\langle g_2,h_0^3, x \rangle$, though
it is also sometimes used when this bracket is not defined.  The operators
$\Delta$ and $\Delta_1$ are given by products with non-permanent
cycles $b_{03}^2$ and $b_{13}^2$, respectively, in the $E_2$-term
of the May spectral sequence.  This is analogous to Tangora's
definition\footnote{Tangora \cite{Tan70}*{pp.~32 and~48} notes that
this is not always equal to the periodicity operator $Px = \langle
h_3, h_0^3, x\rangle$. For example see \cite{Tan70}*{Note~3 on p.~48}.}
$Px = (b_{02})^{2} x$.  They are precursors of Toda brackets
discussed in Section~\ref{sec:operators} in the sense that, when
the brackets are defined, they often compute them.  In fact, each
could be interpreted as two of three distinct brackets which must
sum to zero by the Jacobi identity.  We discuss these brackets in
Section~\ref{sec:operators}.

In a few bidegrees they encounter classes with no name under this
system. They adopt the notation $x_{n,s}$ for such a class if it
is the unique such class in bidegree $(s,s+n)$.   In bidegree
$(10,10+94)$ there are two such classes which they call $x_{94,10}$
and $y_{94,10}$.

\subsection{Cohomological degrees up to $5$}
Recall the theorems of Wang and Palmieri:

\begin{theorem}[\cite{Wan67} and \cite{Pal07}]
For $s<4$, the homomorphism $Sq^0 : \Ext^s \lra \Ext^s$ is injective.
When $s=4$, its kernel is $\langle h_0^4 \rangle$.
\end{theorem}

This makes the ``$Sq^0$-families'' $\{ x, Sq^0(x), (Sq^0)^2(x),
\ldots\}$ especially useful in low cohomological degrees.  

\begin{remark}
Because $h_1^4 = 0$, 
$Sq^0 : \Ext^5 \lra \Ext^5$ must send the nonzero elements
$h_0^5$ and $h_0^4 h_i$, $i \geq 4$, to 0.   In the range we have calculated,
the only other element in its kernel is $Ph_2$, reflecting
$Sq^0(Ph_2) = h_3 g = 0$.
\end{remark}

Chen (\cite{Che11}*{Theorem 1.2 and Theorem 1.3}) gives a complete
description of $\Ext^s$ for $s \leq 5$, building on the work
of Adams~\cite{Ada60}, Wang~\cite{Wan67} and Lin~\cite{Lin08}.

\begin{theorem}[\cite{Che11}*{Theorems 1.2 and 1.3}]
An $\bF_2$-base for the indecomposable elements in cohomological degrees
$s\leq 5$ is as follows. In each, $i$ runs over all $i \geq 0$.
\begin{description}
\item[$\Ext^1$]
$h_i$.
\item[$\Ext^3$]
$c_i$.
\item[$\Ext^4$]
$d_i$, $e_i$, $f_i$, $g_{i+1}$, $p_i$, $D_3(i)$ and $p_i'$.
\item[$\Ext^5$]
$Ph_1$, $Ph_2$, $n_i$, $x_i$, $D_1(i)$, $H_1(i)$, $Q_3(i)$,
$K_i$, $J_i$, $T_i$, $V_i$, $V_i'$ and~$U_i$.
\end{description}
\end{theorem}

Note that $V_0'$ and $U_0$ are in the $252$ and $260$ stems,
respectively, so are beyond the range of our computation.  For each
of the other families in these lists, at least the first member of
the family lies in the range of our calculations.  Chen adopts the
notation $D_1(i)$, et cetera, for the $i^{\mathrm{th}}$ member of
a $Sq^0$-family in order to avoid double subscripts.  We extend
this practice into higher cohomological degrees as noted in the
next section.  In each family except $\{g_1 = g, g_2, \ldots\}$,
the family starts with the $0^{\mathrm{th}}$ element.

Since Chen's list is (mostly) in terms of $Sq^0$-families, we need
only consider the first member of each family in order to establish
the relation between his classes and ours.  Thirteen of the families
start in a bidegree which is $1$-dimensional over $\bF_2$, so that
the correspondence is clear for them.   The remaining families are
$f_i$, $n_i$, $H_1(i)$, $Q_3(i)$ and $T_i$, and we consider each
of these individually.  

\begin{remark}\leavevmode
\label{Chentomin}
\begin{enumerate}
\item
It is long established practice to define $f_0 = Sq^1(c_0)$, because
that allows us to take advantage of $H_\infty$ ring spectrum relations
and differentials.  We do not know whether Chen's definition of
$f_0$ in terms of the Lambda algebra equals this or $f_0 + h_1^3
h_4$.
\item
We choose $n_0$ to be the $h_0$-annihilated class in bidegree
$(5,5+31)$.  According to Chen's unpublished
preprint~\cite{Che12}*{Thm.~1.7} this agrees with his Lambda algebra
definition of $n_0$.

\item
In two remaining cases, $H_1(0)$ and $Q_3(0)$,  their May spectral
sequence definition specifies a coset they must lie in.  Chen's
Lambda algebra classes with these names lie in the specified cosets.
The precise elements in these cosets are determined by the classes
$h_0 H_1(0)$, $h_4 H_1(0)$ and  $h_3 Q_3(0)$.  (That this suffices
can be checked using {\tt all.products}.) These three products are
shown to be zero in Chen's~\cite{Che12}*{Thm.~1.7}.  Thus, the May spectral
sequence indeterminacy for $H_1(0)$ and $Q_3(0)$ is the indeterminacy
reported in Table~\ref{concordance}, but Chen's unpublished results
allow the more precise correspondence $H_1(0) = 5_{32}$, $Q_3(0) =
5_{39}$, $H_1(1) = 5_{81}$ and $Q_3(1) = 5_{91}$.

\item
The remaining case, $T_0$, has only Chen's Lambda algebra definition
from~\cite{Che11}.  The precise indecomposable element in this
bidegree is determined by the values of $h_1T_0$ and $h_4 T_0$.
As above, that these suffice can be checked using {\tt all.products},
and both are shown to be zero in Chen's~\cite{Che12}*{Thm.~1.7}.  
The indeterminacy reported in Table~\ref{concordance} is the set
of decomposables, but Chen's unpublished results
allow the more precise correspondence $T_0 = 5_{93}$.

\end{enumerate}
\end{remark}

\subsection{Cohomological degrees greater than $5$}

In higher cohomological degrees, names for the elements of $\Ext$
come from Tangora's 1970 calculation \cite{Tan70} of the $E_\infty$-term
of the May spectral sequence.  His Appendix 1 lists its indecomposables
in the range $t-s \leq 70$ (omitting classes of the form $P^k a$).
In three cases, known hidden extensions
between the associated graded and $\Ext$ allow us to ignore these
elements. These are $s = h_0 r$, $S_1 \in \langle h_1 x', h_0 R_1
\rangle$ and $g_2' \in \langle h_1 B_{21}, h_0^2 B_{4}(0) \rangle$.

Some elements named by Tangora have $Sq^0$ which is decomposable.
In those cases we keep Tangora's (often unsubscripted) notation.  Others
appear to be the start of $Sq^0$-families, with indecomposable
members in the range of our calculation.  In these cases, we add a
subscript $0$ or suffix $(0)$ to Tangora's name for the class, 
except for the $g_i$ family, which starts with $g_1 = g$.

For unsubscripted elements that start $Sq^0$-families, like Tangora's
$m$ or $A$, we adopt the usual $a_{i+1} = Sq^0(a_i)$ notation for the
subsequent elements of the family.  For subscripted elements like the
$B_i$ or $H_1$ we use Chen's suffix notation, $Z(i+1) = Sq^0(Z(i))$, to
avoid collision with other names used by Tangora.  The suffix notation
is also applied for $C$, $C''$ and~$G$, due to the prior presence of
classes $C_0$ and~$G_0$.

Collisions would occur because subscripts in Tangora's subscripted
capital letter classes do not indicate membership in a $Sq^0$-family.
In particular, the $B_1$ through~$B_5$, the $C$ and $C_0$, et cetera,
are not related
by $Sq^0$.  Nor are $q$ in the $32$-stem and $q_1$ in the $64$-stem
related in this manner.  This is solved by Chen's $Z(i)$ notation.

If there is indeterminacy in the $s_g$ name for the first member
of a $Sq^0$-family, it is generally inherited by the subsequent
members.  $X_2(1)$, $E_1(1)$, $C_0(1)$, $G_{21}(1)$ and $B_4(1)$
are exceptions:  $Sq^0$ annihilates the indeterminacy in the
descriptions of $X_2(0)$, $E_1(0)$, $C_0(0)$, $G_{21}(0)$ and $B_4(0)$.

We have not listed all the indecomposable $s_g$ which can be described
using the Adams periodicity operators, but we have included those
appearing in Tangora's list of indecomposables.  Similarly, we have
not listed all the indecomposable $s_g$ which can be described using
the ``Mahowald operator'' $M(x) = \langle g_2, h_0^3, x \rangle$ discussed
by Margolis--Priddy--Tangora~\cite{MPT70} and Isaksen~\cite{Isa}, but we
have noted that seven classes $B_1$, $B_2$, $B_3$, $B_{21}$, $B_{22}$,
$B_{23}$ and $G_{11}$ in Tangora's list could be so described using the
modified operator $M'(x) = \langle h_0, h_0^2 g_2, x \rangle$.
See Section~\ref{sec:operators} for further discussion
of how to use our data to extend this into the full range of our
calculation.

In a very few cases the indeterminacy reported in Table~\ref{concordance}
is greater than the inherent indeterminacy in a May spectral sequence
definition.  These are the bidegrees $(t-s,s) = (63,7)$,
$(67,9)$ and~$(66,10)$, where there are two indecomposables in the same
bidegree, and possibly $(t-s,s) = (141,5)$, where we do not know a May
spectral sequence definition of the indecomposable class.


\begin{longtable}{|>{$}r<{$}>{$}r<{$}>{$}r<{$}|>{$}l<{$}|>{$}l<{$}|}
\caption{Concordance between indecomposable $s_g$ and other notations.
        Elements $s_g$ for which we do not have a traditional name
	are omitted.
        \label{concordance}} \\
\hline
t-s & s & g & {\mathrm{Tangora, Chen}} & {\mathrm{Note}} \\
\hline
\endfirsthead
\caption{Concordance between indecomposable $s_g$ and other notations.}\\
\hline
t-s & s & g &  {\mathrm{Tangora, Chen}} & {\mathrm{Note}} \\
\hline
\endhead
\hline
\endfoot

0 & 1 & 0 & 	h_0 &  \ref{uniq} \\
1 & 1 & 1 & 	h_1 &  \ref{sq0} \\
3 & 1 & 2 & 	h_2 &  \ref{sq0} \\
7 & 1 & 3 & 	h_3 &  \ref{sq0} \\
15 & 1 & 4 & 	h_4 &  \ref{sq0} \\
31 & 1 & 5 & 	h_5 &  \ref{sq0} \\
63 & 1 & 6 & 	h_6 &  \ref{sq0} \\
127 & 1 & 7 & 	h_7 &  \ref{sq0} \\
\hline
8 & 3 & 3 & 	c_0 &  \ref{uniq} \\
19 & 3 & 9 & 	c_1 &  \ref{sq0} \\
41 & 3 & 19 & c_2 &  \ref{sq0} \\
85 & 3 & 34 & c_3 &  \ref{sq0} \\
173 & 3 & 55 & c_4 &  \ref{sq0} \\
\hline
14 & 4 & 3 & d_0 &  \ref{uniq} \\
17 & 4 & 5 & 	e_0 &  \ref{uniq} \\
18 & 4 & 6 & f_0 & \ref{BNT}  \\
20 & 4 & 8 & g_1 &  \ref{uniq} \\
32 & 4 & 13 & d_1 &  \ref{sq0} \\
33 & 4 & 14 & p_0 &  \ref{uniq} \\
38 & 4 & 16 & e_1 &  \ref{sq0} \\
40 & 4 & 19 & f_1 &  \ref{sq0} \\
44 & 4 & 22 & g_2 &  \ref{sq0} \\
61 & 4 & 26 & D_3(0) &  \ref{uniq} \\
68 & 4 & 32 & d_2 &  \ref{sq0} \\
69 & 4 & 33 & p_0' &  \ref{uniq} \\
70 & 4 & 34 & p_1 &  \ref{sq0} \\
80 & 4 & 40 & e_2 &  \ref{sq0} \\
84 & 4 & 44 & f_2 &  \ref{sq0} \\
92 & 4 & 48 & g_3 &  \ref{sq0} \\
126 & 4 & 53 & D_3(1) &  \ref{sq0} \\
140 & 4 & 65 & d_3 &  \ref{sq0} \\
142 & 4 & 67 & p_1' &  \ref{sq0} \\
144 & 4 & 69 & p_2 &  \ref{sq0} \\
164 & 4 & 79 & e_3 &  \ref{sq0} \\
172 & 4 & 84 & f_3 &  \ref{sq0} \\
188 & 4 & 89 & g_4 &  \ref{sq0} \\
\hline
9 & 5 & 1 & Ph_1 &  \ref{uniq} \\
11 & 5 & 2 & Ph_2 &  \ref{uniq} \\
31 & 5 & 13 & n_0 & \ref{h0}  \\
37 & 5 & 17 & x_0 &  \ref{uniq} \\
52 & 5 & 30 & D_1(0) &  \ref{uniq} \\
62 & 5 & 32(33,34) & H_1(0) & \ref{H1}  \\
67 & 5 & 38 & n_1 &  \ref{sq0} \\
67 & 5 & 39(38) & Q_3(0) & \ref{Q3}  \\
79 & 5 & 50 & x_1 &  \ref{sq0} \\
109 & 5 & 75 & D_1(1) &  \ref{sq0} \\
125 & 5 & 77 & K_0 &  \ref{uniq} \\
128 & 5 & 80 & J_0 &  \ref{uniq} \\
129 & 5 & 81(82,83) & H_1(1) &  \ref{sq0} \\
139 & 5 & 90 & n_2 &  \ref{sq0} \\
139 & 5 & 91(90) & Q_3(1) &  \ref{sq0} \\
141 & 5 & 93(94,95) & T_0 & \ref{T0}  \\
156 & 5 & 108 & V_0 &  \ref{uniq} \\
163 & 5 & 115 & x_2 &  \ref{sq0} \\
\hline
30 & 6 & 10 & r_0 &  \ref{uniq} \\
32 & 6 & 12 & q &  \ref{uniq} \\
36 & 6 & 14 & t_0 &  \ref{uniq} \\
38 & 6 & 16 & y &  \ref{BNT} \\
50 & 6 & 27 & C(0) &  \ref{uniq} \\
54 & 6 & 30 & G(0) &  \ref{uniq} \\
58 & 6 & 31 & D_2 &  \ref{uniq} \\
61 & 6 & 32+33 & A_0' &  \ref{AAprime} \\
61 & 6 & 33 & A_0 &  \ref{AAprime} \\
64 & 6 & 38 & A_0'' &  \ref{uniq} \\
66 & 6 & 40 & r_1 &  \ref{sq0} \\
78 & 6 & 56 & t_1 &  \ref{sq0} \\
106 & 6 & 87 & C(1) &  \ref{sq0} \\
114 & 6 & 92 & G(1) &  \ref{sq0} \\
128 & 6 & 102 & A_1+h_2K_0 = A_1'+h_0J_0 & \ref{6102} \\
134 & 6 & 110 & A_1'' &  \ref{sq0} \\
138 & 6 & 115 & r_2 &  \ref{sq0} \\
162 & 6 & 156 & t_2 &  \ref{sq0} \\
\hline
16 & 7 & 3 & Pc_0 &  \ref{uniq} \\
23 & 7 & 5 & i &  \ref{uniq} \\
26 & 7 & 6 & j &  \ref{uniq} \\
29 & 7 & 7 & k &  \ref{uniq} \\
32 & 7 & 10 & \ell &  \ref{uniq} \\
35 & 7 & 12 & m_0 &  \ref{uniq} \\
46 & 7 & 20 & B_1(0) = M'h_1 &  \ref{uniq} \\
48 & 7 & 22(23) & B_2(0) = M'h_2 &  \ref{B2} \\
57 & 7 & 27 & Q_2(0)  & \ref{uniq} \\ 
60 & 7 & 29 & B_3 = M'h_4  & \ref{uniq} \\ 
63 & 7 & 33(35) & X_2(0)  & \ref{CprimeX2} \\ 
63 & 7 & 34(33,35) & C'  & \ref{CprimeX2} \\ 
66 & 7 & 40(41) & G_0(0) & \ref{G0}  \\ 
77 & 7 & 56 & m_1 + h_1 6_{53} &  \ref{sq0} \\
99 & 7 & 85 & B_1(1) &  \ref{sq0} \\
103 & 7 & 90(91) & B_2(1) &  \ref{sq0} \\
121 & 7 & 101 & Q_2(1) +h_6 D_2 & \ref{sq0} \\  
133 & 7 & 124 & X_2(1) + h_3 6_{97} + h_1 6_{107} & \ref{modmsq} \\
133 & 7 & 125 & X_2(1) + h_1 6_{107} & \ref{modmsq} \\
139 & 7 & 137(138) & G_0(1) & \ref{sq0} \\ 
161 & 7 & 184 & m_2 + h_2 6_{149} + h_2 h_7 n &  \ref{sq0} \\
\hline
22 & 8 & 3 & Pd_0 &  \ref{uniq} \\
25 & 8 & 5 & Pe_0 &  \ref{uniq} \\
46 & 8 & 20 & N & \ref{uniq} \\
62 & 8 & 32+33(34) &  E_1(0)  & \ref{C0E1}  \\ 
62 & 8 & 33(34) &  C_0(0)  & \ref{C0E1}  \\ 
68 & 8 & 43(44) & G_{21}(0) & \ref{G21} \\
69 & 8 & 46 & PD_3(0) & \ref{Ptxt} \\
132 & 8 & 139+140 & C_0(1) + h_2^2 6_{97} & \ref{sq0} \\
132 & 8 & 140 & E_1(1) & \ref{sq0} \\
144 & 8 & 176 & G_{21}(1) + h_1^3 h_7 d_0 & \ref{sq0} \\
\hline
17 & 9 & 1 & P^2h_1 &  \ref{uniq} \\
19 & 9 & 2 & P^2h_2 &  \ref{uniq} \\
39 & 9 & 18 & u &  \ref{uniq} \\
42 & 9 & 19 & v &  \ref{uniq} \\
45 & 9 & 20 & w &  \ref{uniq} \\
60 & 9 & 29(30) & B_4(0) & \ref{B4} \\
61 & 9 & 31 & X_1 & \ref{uniq} \\
67 & 9 & 39(40) & C''(0) & \ref{CprimeprimeX3} \\
67 & 9 & 40(39) & X_3 & \ref{CprimeprimeX3} \\
67 & 9 & 39,40 & C''(0), X_3 & \ref{CprimeprimeX3} \\
129 & 9 & 145 & B_4(1) + h_2 8_{118} &  \ref{sq0} \\
143 & 9 & 197(199+200) & C''(1) &  \ref{sq0} \\
\hline
41 & 10 & 14 & z & \ref{uniq}  \\
53 & 10 & 18 & x' & \ref{uniq}  \\
54 & 10 & 19(20) & R_1 & \ref{R1}   \\
56 & 10 & 22(21) & Q_1  & \ref{Q1}   \\
59 & 10 & 24 & B_{21} = M'd_0 & \ref{uniq}  \\
62 & 10 & 27(28,29) & R & \ref{B22R}  \\  
62 & 10 & 28(29) & B_{22} \ni M'e_0 & \ref{B22R}  \\
64 & 10 & 32(33) & q_1 & \ref{q1}  \\
65 & 10 & 34 & B_{23}(0) = M'g_1 & \ref{uniq} \\  
66 & 10 & 35 + 36 & B_5  & \ref{D2primeB5} \\
66 & 10 & 35 {\ \text{or}\ } 36 & D_2' & \ref{D2primeB5} \\
66 & 10 & 36 & PD_2 & \ref{Ptxt} \\
69 & 10 & 40 & PA_0  & \ref{uniq} \\
140 & 10 & 196+197 & B_{23}(1)+ h_1 9_{178} & \ref{B23}  \\
\hline
24 & 11 & 3 & P^2c_0 &  \ref{uniq} \\
34 & 11 & 7 & Pj &  \ref{uniq} \\
67 & 11 & 35(36) & C_{11} & \ref{C11} \\
\hline
30 & 12 & 3 & P^2d_0 &  \ref{uniq} \\
33 & 12 & 5 & P^2e_0 &  \ref{uniq} \\
\hline
25 & 13 & 1 & P^3h_1 &  \ref{uniq} \\
27 & 13 & 2 & P^3h_2 &  \ref{uniq} \\
47 & 13 & 14 & Q  & \ref{PuQ}  \\
47 & 13 & 14+15 & Q' = Q + Pu  & \ref{PuQ}  \\
47 & 13 & 15 & Pu  & \ref{PuQ} \\
50 & 13 & 16 & Pv  & \ref{uniq} \\
65 & 13 & 29(28) & R_2 & \ref{R2} \\
68 & 13 & 30(31) & G_{11} = M'i & \ref{G11}  \\
69 & 13 & 32 & W_1 & \ref{uniq}  \\
\hline
64 & 14 & 23 & PQ_1 & \ref{Ptxt} \\
\hline
70 & 17 & 26(27) & R_1' & \ref{R1prime} \\
\hline
69 & 18 & 20 & \text{$P^2x'$ (ill-defined)} & \ref{P2xprime} \\
\hline
8k+1 & 4k+1 & 1 & P^kh_1 &  \ref{per} \\
8k+3 & 4k+1 & 2 & P^kh_2 &  \ref{per} \\
\hline
8k+8 & 4k+3 & 3 & P^kc_0 &  \ref{per} \\
\hline
8k+14 & 4k+4 & 3 & P^kd_0 &  \ref{per} \\
8k+17 & 4k+4 & 5 & P^ke_0 &  \ref{per} \\
\hline
8k+23 & 4k+7 & 5 & P^ki \text{\ (for $k$ even)} & \ref{per} \\
8k+26 & 4k+7 & 6 & P^kj \text{\ (for $k$ even)} & \ref{per} \\
8k+26 & 4k+7 & 7 & P^kj \text{\ (for $k$ odd)} & \ref{per} \\
\hline
8k+39 & 4k+9 & 18 & P^ku \text{\ (for $k$ even)} & \ref{per} \\
8k+39 & 4k+9 & 15 & P^ku \text{\ (for $k$ odd)} & \ref{per} \\
8k+42 & 4k+9 & 19 & P^kv \text{\ (for $k$ even)} & \ref{per} \\
8k+42 & 4k+9 & 16 & P^kv \text{\ (for $k \equiv 1 \mod 4$)} & \ref{per} \\
8k+42 & 4k+9 & 17 & P^kv \text{\ (for $k \equiv 3 \mod 4$)} & \ref{per} \\
\end{longtable}

Notes:

\begin{enumerate}
\item
\label{uniq}
There is a unique nonzero element in this bidegree.
\item
\label{sq0}
This is $Sq^0$ of an element we have already identified.
\item
\label{BNT}
We choose to let $f_0 = Sq^1(c_0)$ and $y = Sq^2(f_0)$, which are shown
in~\cite{BNT} to be $4_6$ and $6_{16}$, respectively.
\item
\label{h0}
This is the unique nonzero element in this bidegree
whose $h_0$ multiple is $0$.
\item
\label{H1}
The May spectral sequence definition of $H_1 = H_1(0)$ has indeterminacy
spanned by $5_{33} = h_1D_3$ and $5_{34} = h_0^3 h_5^2$, so that
$H_1(0)$ must be the indecomposable $5_{32}$ modulo them.  The
Lambda algebra class which Chen defines as $H_1(0)$ in~\cite{Che11}
satisfies $h_0 H_1(0) = 0$ and $h_4 H_1(0) = 0$, according to Chen's
preprint~\cite{Che12}*{Thm.~1.7},  eliminating the possible summands
$5_{34}$ and $5_{33}$, respectively.
\item
\label{Q3}
The May spectral sequence definition of $Q_3 = Q_3(0)$ has indeterminacy
spanned by $5_{38} = n_1$, so that $Q_3(0)$ must be $5_{39}$ modulo
$5_{38}$.  The Lambda algebra class which Chen defines as $Q_3(0)$
in~\cite{Che11} satisfies $h_3 Q_3(0) = 0$, according to Chen's
preprint~\cite{Che12}*{Thm.~1.7},  eliminating the possible summand
$5_{38}$.
\item
\label{T0}
The indecomposable $T_0$ must be the indecomposable $5_{93}$ modulo
the decomposables $5_{94} = h_7 d_0$ and $5_{95} = h_1 d_3$.  If
$T_0 = 5_{93} + \alpha 5_{94} + \beta 5_{95}$, then $h_1 T_0 =
\alpha 6_{126}$ while $h_4 T_0 = \beta 6_{145}$.  These products
are both zero according to Chen's preprint~\cite{Che12}*{Thm.~1.7},
so that $T_0 = 5_{93}$.
\item
\label{AAprime}
Tangora \cite{Tan70} shows $h_0 A = h_2 D_2$, hence $A = 6_{33}$.
He also shows $h_0^2 A' = 0$, so that $A' = 6_{32} + 6_{33}$.  We
write them as $A_0$ and $A_0'$, since $Sq^0$ is nonzero on both.
\item
\label{6102}
$A_1 = Sq^0(A_0) = 6_{102} + 6_{103}$ and 
$A_1' = Sq^0(A_0') = 6_{102} + 6_{104}$, while
$6_{103} = h_2 K_0$ and $6_{104} = h_0 J_0$.
\item
\label{B2}
The May spectral sequence definition of $B_2$ has indeterminacy
$7_{23} = h_0^2h_5e_0$, so that $B_2 \in \{7_{22}, 7_{22}+7_{23}\}$.
The value of $M'h_2$ is exactly the same set. 
\item
\label{CprimeX2}
Bidegree $(7,7+63)$ is spanned by the two indecomposables $7_{33}$
and $7_{34}$ together with $7_{35} = h_0^6h_6$.  We have $h_2 7_{33}
= 0$ and $h_2 7_{34} = 8_{41}$.  In the May spectral sequence there are
indecomposables $C'$ and $X_2$.  Tangora~\cite{Tan70} reports that $h_2 C'
= h_0 G_0$ is nonzero, and does not list $h_2 X_2$ as a nonzero value.
Granting that $h_2 X_2 = 0$, it follows that $C' \equiv 7_{34}$ modulo
$\langle 7_{33}, 7_{35} \rangle$ and $X_2 \equiv 7_{33}$ modulo $7_{35}$.
Since $Sq^0(7_{33})$ is indecomposable, we write $X_2 = X_2(0)$.  We do
not know whether $Sq^0(C')$ is indecomposable; if it is we should set $C'
= C'(0)$ and note that $C'(1) \equiv X_2(1)$ modulo decomposables.
\item
\label{G0}
Bidegree $(7,7+66)$ is spanned by the indecomposable $7_{40}$ and
$7_{41} = h_0r_1$.  Hence,  the indecomposable $G_0 = G_0(0)$ must be $7_{40}$
modulo $7_{41}$.
\item
\label{modmsq}
Since $Sq^0(7_{34}) = 7_{124} + 7_{125} = h_3 6_{97}$, we see that
$7_{124}$ and $7_{125}$ are congruent modulo decomposables, but are
each indecomposable.  Then $Sq^0(7_{33}) = 7_{125} + h_1 6_{107}$
shows that $7_{125} = Sq^0(7_{33}) + h_1 6_{107}$ and
that $7_{124} = Sq^0(7_{33}) + h_1 6_{107} + h_3 6_{97}$.
\item
\label{C0E1}
Bidegree $(8,8+62)$ is spanned by the two indecomposables $8_{32}$ and
$8_{33}$ together with $8_{34} = h_0^6h_5^2$.  Tangora~\cite{Tan70} lists
indecomposables $C_0$ and $E_1$, together with $h_0^6h_5^2$.  He reports
$h_1 E_1 \neq 0$ and (implicitly) $h_1 C_0 = 0$.  He also reports $h_2
C_0 \ne 0$ and (implicitly) $h_2 E_1 \equiv 0$ modulo $h_0^2 h_3 D_2$.
This implies $E_1 \equiv 8_{32} + 8_{33}$ modulo~$8_{34}$ and $C_0
\equiv 8_{33}$ modulo~$8_{34}$.  We set $E_1(0) = E_1$ and
$C_0(0) = C_0$ since $Sq^0(E_1) = 8_{140}$ and $Sq^0(C_0)
= 8_{139} + 8_{140} + 8_{141} + 8_{142}$ are indecomposable.
\item
\label{G21}
Bidegree $(8,8+68)$ is spanned by the indecomposable $8_{43}$ and
the decomposable $8_{44} = h_0h_3 A_0'$.  Hence,  the indecomposable
$G_{21}$ must be $8_{43}$ modulo $8_{44}$.
\item
\label{Ptxt}
This is the unique element in the bracket
$\langle h_3, h_0^4, - \rangle$.
\item
\label{B4}
The May spectral sequence definition of $B_4 = B_4(0)$ has indeterminacy
spanned by $9_{30} = h_0^2 B_3$, so that $B_4(0)$ must be the
indecomposable $9_{29}$ modulo $9_{30}$
\item
\label{CprimeprimeX3}
Bidegree $(9,9+67)$ is spanned by the two indecomposables $9_{39}$
and~$9_{40}$.  We have $h_0 9_{39} = 0$, $h_0 9_{40} = 10_{37}$, $h_2
9_{39} = 10_{41}$ and $h_2 9_{40} = 0$.  In the May spectral sequence,
there are indecomposables $C''$ and $X_3$.  The relations reported
by Tangora~\cite{Tan70} have $h_0 X_3 \ne 0$ and $h_2 C'' \ne 0$.
It follows that $C'' \equiv 9_{39}$ modulo~$9_{40}$ and $X_3 \equiv
9_{40}$ modulo~$9_{39}$, but, of course, $C'' \ne X_3$.
Since $Sq^0(C'')$ is indecomposable, we write $C'' = C''(0)$.  We do
not know whether $Sq^0(X_3)$ is indecomposable; if it is we should set $X_3
= X_3(0)$ and note that $C''(1) \equiv X_3(1)$ modulo decomposables.
\item
\label{R1}
Bidegree $(10,10+54)$ is spanned by the indecomposable $10_{19}$
and the decomposable $10_{20} = h_0^2 h_5 i$, while Tangora's calculation
has this bidegree spanned by $h_0^2 h_5 i$ and $R_1$.  Hence, $R_1$ must
be $10_{19}$ modulo $10_{20}$.
\item
\label{Q1}
Bidegree $(10,10+56)$ is spanned by the indecomposable $10_{22}$
and the decomposable $10_{21} = g_1t_0$, while Tangora's calculation
has this bidegree spanned by $gt$ and $Q_1$.  Hence, $Q_{1}$ must
be $10_{22}$ modulo $10_{21}$.
\item
\label{B22R}
Bidegree $(10,10+62)$ is spanned by the decomposable $10_{29} = h_1
X_1 = PG$ together with the indecomposables $10_{27}$ and $10_{28}$.
In the May spectral sequence it is spanned by $R$, $B_{22}$ and $PG$,
in order of May filtration, so that $B_{22}$ is defined modulo~$PG$, while
$R$ is only defined modulo the other two.  The relation $h_0 B_{22}
= d_0 B_2$ holds in the May spectral
sequence by~\cite{Tan70}.  Since $h_0 10_{28} = d_0 B_2$, while
$h_0 10_{27}$ is not divisible by $d_0$ in $\Ext$, $B_{22}$ must
be $10_{28}$ modulo~$10_{29}$.
Since $R$ is linearly independent of $B_{22}$ and $PG$,
it must be $10_{27}$ modulo $\langle 10_{28}, 10_{29} \rangle$.
\item
\label{q1}
Tangora's $q_1$ in (10,10+64) is not $Sq^0(q)$, which is instead the
decomposable $6_{46} = Sq^0(6_{12}) = h_2 Q_3(0)$.  The May spectral
sequence definition of $q_1$ has indeterminacy $10_{33} = h_0^2h_3Q_2(0) =
h_1^2E_1(0)$ so that $q_1$ is $10_{32}$ modulo $10_{33}$.  
\item
\label{D2primeB5}
Bidegree $(10,10+66)$ is spanned by the two indecomposables $10_{35}$
and $10_{36}$.  We have $h_0 (10_{35} + 10_{36})= h_1 B_{23}(0) =
h_2^2 B_4(0)$, while $h_0^2 10_{35} = h_0^2 10_{36} = 12_{32} =
h_1^2 q_1$.  In the May spectral sequence, there are indecomposables
$D_2'$ and $B_5$.  The relations reported by Tangora~\cite{Tan70}
include $h_0 B_5 = h_1 B_{23} = h_2^2 B_4$, which require $B_5 =
10_{35} + 10_{36}$.   We then have $D_2' = 10_{35}$ or $10_{36}$.
In any case $\langle D_2',B_5\rangle = \langle 10_{35}, 10_{36} \rangle$.
\item
\label{B23}
$B_{23}(1) = Sq^0(B_{23}(0)) = 10_{196} + 10_{197} + 10_{199}$ and
$10_{199} = h_1 9_{178} = h_1 9_{179}$.   Both $10_{196}$ and $10_{197}$
are indecomposable.
\item
\label{C11}
Bidegree $(11,11+67)$ is spanned by the indecomposable $11_{35}$
and $11_{36} = h_0^2 X_3 = h_0h_3B_4(0) = ig_2 = r_0x_0$. Hence,
the indecomposable $C_{11}$ must be $11_{35}$ modulo $11_{36}$.
\item
\label{PuQ}
Bidegree $(13,13+47)$ is spanned by the indecomposables $13_{14}$
and $13_{15}$, while Tangora's calculation has this bidegree spanned
by $Pu$ and $Q$ with $h_1 Q \neq 0$.  The brackets file shows that
$Pu = \langle h_3, h_0^4, u \rangle = 13_{15}$.  Since $h_1 13_{14}
= h_1 13_{15}$, we must have $Q = 13_{14}$.  Tangora defines $Q' = Q + Pu$.
\item
\label{R2}
Bidegree $(13,13+65)$ is spanned by the indecomposable $13_{29}$
and the decomposable $13_{28} = g_1 w = r_0 m $, while Tangora's calculation
has this bidegree spanned by $gw$ and $R_2$.  Hence, $R_2$ must
be $13_{29}$ modulo $13_{28}$.
\item
\label{G11}
Bidegree $(13,13+68)$ is spanned by the indecomposable $13_{30}$
and the (highly) decomposable $13_{31} = h_0^5 G_{21} = h_0h_5d_0i$. 
Hence, the indecomposable $G_{11}$ must be $13_{30}$ modulo $13_{31}$.
This coset is $M'i$ since $13_{31}$ is in the indeterminacy of $M'i$.
\item
\label{R1prime}
Bidegree $(17,17+65)$ is spanned by the indecomposable $17_{26}$
and the decomposable $17_{27} = d_0^2 v $, while Tangora's calculation
has this bidegree spanned by $Pe_0w = d_0^2 v$ and $R_1'$.  Hence, $R_1'$ must
be $17_{26}$ modulo $17_{27}$.
\item
\label{P2xprime}
Bidegree $(18,18+69)$ is spanned by the indecomposable $18_{20}$ together
with $18_{21} = P(gz)$.  Tangora writes $P^2x' = \langle h_4, h_0^8, x'
\rangle$ for an indecomposable in this bidegree, but this is not defined
as a Toda bracket, since $h_0^8 x' \ne 0$.
\item
\label{per}
By Adams periodicity $P^kh_1 = (4k+1)_1$, $P^kh_2 = (4k+1)_2$, $P^kc_0 =
(4k+3)_3$, $P^kd_0 = (4k+4)_3$ and $P^k e_0 = (4k+4)_5$.  Likewise,
$P^ki = (4k+7)_5$, $P^kj = (4k+7)_6$, $P^ku = (4k+9)_{18}$ and~$P^kv =
(4k+9)_{19}$ for $k$ even, $P^kj = (4k+7)_7$ and~$P^ku = (4k+9)_{15}$
for $k$ odd, $P^kv = (4k+9)_{16}$ for $k \equiv 1 \mod 4$ and $P^kv =
(4k+9)_{17}$ for $k \equiv 3 \mod 4$.
\end{enumerate}

\section{Operators}
\label{sec:operators}

In this section, we describe the files {\tt P.txt}, {\tt P2.txt},
{\tt P4.txt} and {\tt MM.txt}, in which we collect the data needed to compute
the Adams periodicity operators~$P^k$ for~$k=1, 2, 4$,
and the Mahowald operator~$M$, then make some general remarks
about such operators.

\subsection{Adams and Mahowald operators}

These are defined by
\begin{enumerate}
\item
$Px = \langle h_3, h_0^4, x \rangle$,
\item
$P^2x = \langle h_4, h_0^8, x \rangle$,
\item
$P^4x = \langle h_5, h_0^{16}, x \rangle$ and
\item
$Mx = \langle g_2, h_0^{3}, x \rangle$.
\end{enumerate}
The Mahowald operator $Mx = \langle g_2, h_0^3, x \rangle$
is not of the form that it is immediately evident in the {\tt
brackets.sym} files, but its variant $M'(x) = \langle h_0, h_0^2
g_2, x \rangle$ is.  Both contain $\langle  h_0 g_2, h_0^2, x \rangle
$, when defined.  

\begin{landscape}
\begin{figure}
\vskip 1.2cm 
\includegraphics[scale=0.20]{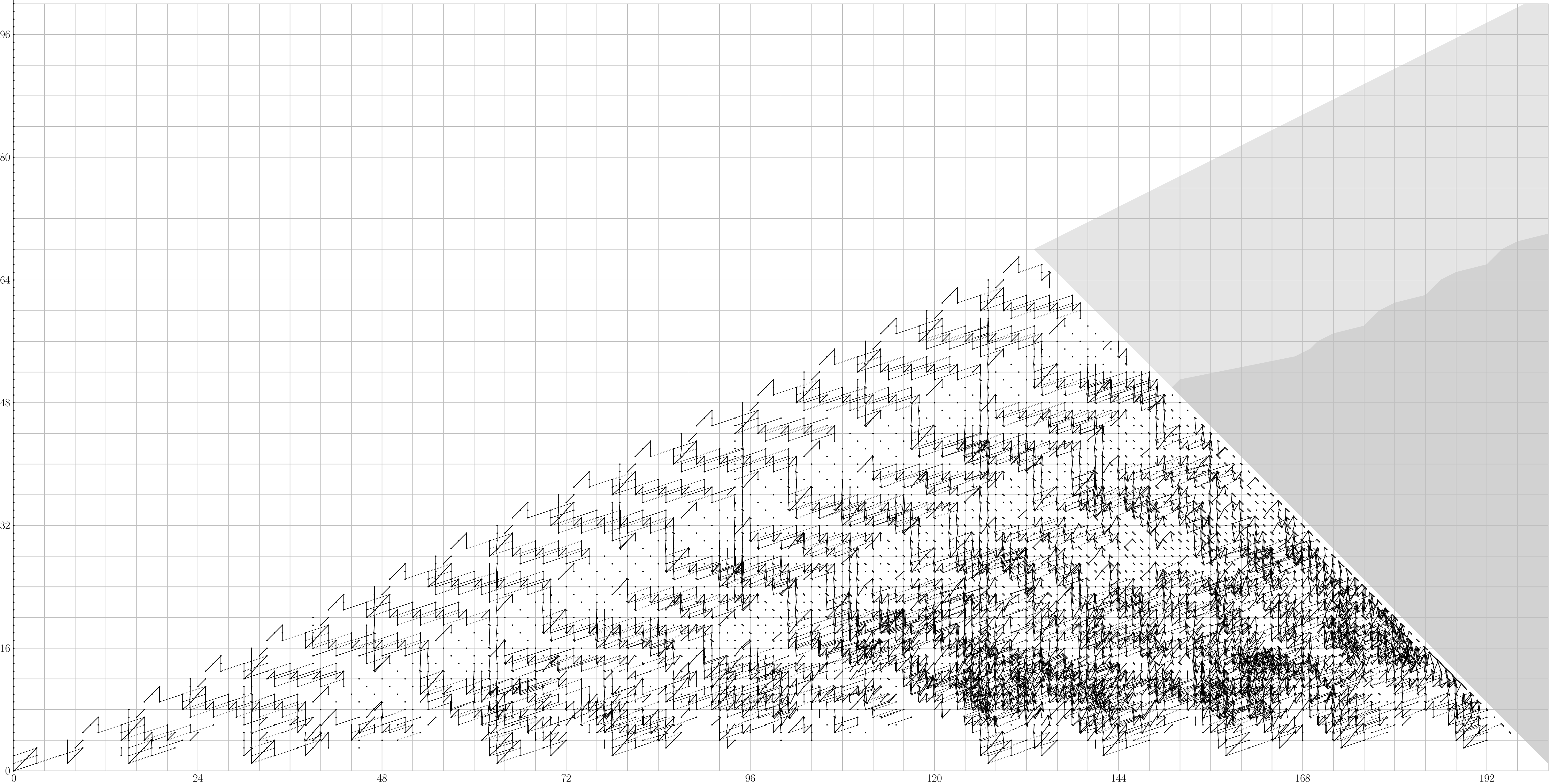}
\caption{$\Ext_A^{s,t}(\bF_2, \bF_2)$ for $t \le 200$}
\label{Ext-A-F2-F2-0-200+gray}
\end{figure}
\end{landscape}

In Figure~\ref{Ext-A-F2-F2-0-200+gray}, the lighter gray region
with $t>200$ and $t-s\le200$ shows a range of bidegrees in which
$\Ext_{\cA}^{s,t}(\bF_2, \bF_2)$ is determined by Adams periodicity and
our calculations for $t\le200$.  The darker gray region shows bidegrees
where further calculation would be required to determine these groups.

Recall from Proposition~\ref{prop:brackets}
that, if defined, the bracket $\langle h_i, (s_0)_{g_0} , (s_1)_{g_1}
\rangle$ is the sum of those $s_g$ with $s = s_0 + s_1$
such that the file {\tt s1\_g1/brackets.sym} contains a line
\begin{verbatim}
s_g in < hi, g0, s1_g1 >
\end{verbatim}

Thus, the values of the operators $P$, $P^2$, $P^4$ and $M'$ are recognized
by the presence of lines
\begin{enumerate}
\item
\begin{verbatim}
s_g in < h3, 0, s1_g1 >
\end{verbatim}
\item
\begin{verbatim}
s_g in < h4, 0, s1_g1 >
\end{verbatim}
\item
\begin{verbatim}
s_g in < h5, 0, s1_g1 >
\end{verbatim}
\item
\begin{verbatim}
s_g in < h0, 21, s1_g1 >
\end{verbatim}
\end{enumerate}
with $s=s_1+4$, $s=s_1+8$, $s=s_1+16$ or $s=s_1+6$, respectively.   
For $M'$, we note that  $6_{21} = h_0^2 g_2$.   

The files {\tt P.txt}, {\tt P2.txt}, {\tt P4.txt} and {\tt MM.txt}
collect this information from all the map files, together with 
information about products needed in determining the domain of
definition and the indeterminacy.    
Each file starts with a short header describing the operator 
and the file's organization, then has three sections: 
\begin{description}
\item[(a)] values of the brackets,
\item[(b)] nonzero products which obstruct existence of the bracket, and
\item[(c)] nonzero products which give the indeterminacy.
\end{description}
In general, the indeterminacy in a bracket $\langle a,b,c \rangle$
is $a(\Ext) + (\Ext)c$.   However, the brackets $P$, $P^2$, $P^4$
and $M'$ have indeterminacy $a(\Ext)$.    For the periodicity
operators, this is because $\Ext^{4,12}$, $\Ext^{8,24}$ and
$\Ext^{16,48}$ are zero, so that $(\Ext)c = 0$ in the relevant
bidegree.   For $M' = \langle h_0, h_0^2 g_2, - \rangle$, it is
because $\Ext^{6,51} = \langle h_0 h_5 d_0 \rangle$, so that $(\Ext)c$
is contained in $a(\Ext) = h_0(\Ext)$ in the relevant bidegree.

For example, the file {\tt P.txt} starts
\begin{verbatim}
% Adams operator P, in the range t\le200
%
% (a) Brackets Px = < h3, h0^4, x >.
% (b) Nonzero products h0^4 * x, obstructing existence.
% (c) Nonzero products h3 * y, giving indeterminacy.

% (a) Brackets Px = < h3, h0^4, x >.

5_1 in < h3, 0, 1_1 >
5_2 in < h3, 0, 1_2 >
5_5 in < h3, 0, 1_3 >

6_1 in < h3, 0, 2_1 >
\end{verbatim}

From the first three lines in part (a) we see that, if they
are defined, $Ph_1 = 5_1$, $Ph_2 = 5_2$ and $Ph_3 = 5_5 = h_0^4 h_4$, 
modulo their indeterminacy.  To determine whether they are defined,
we look at part (b), which starts
\begin{verbatim}
% (b) Nonzero products h0^4 * x, obstructing existence.
  ...
  5    0  (  4    0     F2)  1_0
  5    5  (  4    0     F2)  1_4
  5   14  (  4    0     F2)  1_5
  5   35  (  4    0     F2)  1_6
  5   79  (  4    0     F2)  1_7
  ...
\end{verbatim}
We see that there are no lines ending in {\tt 1\_1},
{\tt 1\_2} or {\tt 1\_3}, so that $h_0^4$ annihilates 
all three of these.   Thus $Ph_1$, $Ph_2$ and $Ph_3$ are
all defined, but $Ph_0$, $Ph_4$, \ldots, are not.

Next, we consider the indeterminacy, which equals the
$h_3$ multiples in the bidegree of the bracket.  The bidegrees
in question are $\Ext^{5,5+9}= \langle 5_1 \>$, 
$\Ext^{5,5+11}= \langle 5_2 \>$  and
$\Ext^{5,5+15}= \langle 5_4, 5_5 \>$.  
We look at part (c) to see the $h_3$ multiples among these.
We find
\begin{verbatim}
% (c) Nonzero products h3 * y, giving indeterminacy.
  ...
  4   85  (  1    3     F2)  3_54
  5   11  (  1    3     F2)  4_5
  5   18  (  1    3     F2)  4_12
  5   19  (  1    3     F2)  4_13
  ...
\end{verbatim}
Since we do not find $5_1$, $5_2$, $5_4$ or $5_5$ among the $h_3$-multiples,
the indeterminacy is~0.

For an example of a possible bracket which is, in fact, undefined, consider the entry
\begin{verbatim}
8_91 in < h3, 0, 4_48 >
\end{verbatim}
in part (a).  This does not mean that $P(4_{48}) = 8_{91}$ because
in part (b) we find that $h_0^4 \cdot 4_{48} = 8_{77} \neq 0$:
\begin{verbatim}
  8   77  (  4    0     F2)  4_48
\end{verbatim}

Finally, for an example with nontrivial indeterminacy,
in part (a) we find an entry
\begin{verbatim}
6_5 in < h3, 0, 2_5 >
\end{verbatim}
so that $6_5 \in P(2_5)$, but in part (c) we find 
\begin{verbatim}
  6    5  (  1    3     F2)  5_1
\end{verbatim}
showing that $6_5 = h_3 \cdot 5_1$ is in the indeterminacy.   Hence
$P(2_5) = \{ 0, 6_5\}$.  (We have $6_5 = h_1^2 d_0 = h_3 Ph_1 = c_0^2$.)

\subsection{General remarks}

We finish this section with some observations about operators like
those we have just considered.  Defining $Px = (b_{02})^2x$ in the
May spectral sequence is justified by the differential $d_4(b_{02}^2)
= h_0^4 h_3$ (\cite{Tan70}*{Prop.~4.3}).  However, this accounts
for only part of the definition of a Massey product or Toda
bracket.
It is simplest\footnote{An idealistic treatment would instead consider Toda brackets of
chain maps, or Massey products in $\End(C_*)$, which is homotopy
commutative but not commutative.}
to discuss this in terms of Massey
products in a commutative differential graded algebra over $\bF_2$
like the $E_r$ terms of the May spectral sequence.

Consider classes $a$, $b$ and $x$ satisfying $ab = bx = xa = 0$.
The Jacobi identity says that 
\[
0 \in
\langle a,b,x \rangle +
\langle b,x,a \rangle +
\langle x,a,b \rangle \,.
\]
If we choose $A$, $U$ and $V$ such that
$d(A) = ab$, $d(U) = bx$ and $d(V) = xa = ax$, then we have
\begin{enumerate}
\item
$Ax + aU \in \langle a, b, x \rangle$,
\item
$Ua + bV \in \langle b, x, a \rangle$ and
\item
$Ax + bV \in \langle b, a, x \rangle = \langle x, a, b \rangle$.
\end{enumerate}
Approximating the bracket $\langle a,b,x \rangle$ by $Ax$, in those cases where
$Ax$ is a cycle, fails to distinguish between
$\langle a,b,x \rangle$ and $\langle b,a,x \rangle = \langle x,a,b \rangle$.
These differ by $\langle b,x,a \rangle = \langle a,x,b \rangle$.

This can lead to greater indeterminacy and to anomalies like Tangora's
observation \cite{Tan70}*{Note~3, p.~48} that $h_5 i$ is annihilated  by
$h_0^3$, but $P(h_5 i)$, if defined to be $(b_{02})^2 h_5 i$, has
$h_0^9 \cdot (b_{02})^2 h_5 i \neq 0$.  In fact, consulting {\tt P.txt}
we see that $P(h_5 i) = P(8_{26}) = 0$ with zero indeterminacy.


By using only the precisely defined brackets we limit the indeterminacy
and get the advantages of their good formal behavior.

Finally, let us point out two other operators which may be of use.
They are
\begin{enumerate}
\item
the complex Bott periodicity operator 
$v_1(x) = \langle h_0, h_1, x \rangle$, which acts on $h_1$-annihilated
classes such as the unit in $\Ext_{E(Q_0,Q_1)}(\bF_2,\bF_2) \Lra \pi_*ku$, and
\item
the mod 2 complex Bott periodicity operator 
$v_1'(x) = \langle h_1, h_0, x \rangle$, which acts on $h_0$-annihilated
classes.
\end{enumerate}
In the universal example, the Adams spectral sequence for  
$\pi_*(S \cup_2 e^1 \cup_\eta e^2)$, 
$v_1(0_0)$ is an $h_1$-annihilated class which supports an 
infinite $h_0$-tower, while
$v_1'(0_0)$ is an $h_0$-annihilated class which supports 
$h_1^2$-multiplication.
These are visible in the {\tt brackets.sym} files in the form
\begin{verbatim}
s_g in < h0, 1, s1_g1 >
\end{verbatim}
and
\begin{verbatim}
s_g in < h1, 0, s1_g1 >
\end{verbatim}
respectively.

\section{Validity}

Several checks have been run to test the validity of the data.   
\begin{enumerate}
\item
After the resolution was computed, a separate computation was done
to check that $d^2 = 0$.

\item
A check of the $\bF_2$-dimension of the kernel and image at each
step, to ensure that the image at each step has the same dimension
as the kernel at the previous step. This checks exactness when
combined with the check that $d^2=0$.

\item
A check that the {\tt Map} files are complete.   If this were not
true, an element $s_g$ missing from {\tt s1\_g1/Map} might not be
reported in {\tt all.products} as a term in a product
$(s_0)_{g_0}\cdot(s_1)_{g_1}$ even when it belongs there.

\item
A check that the {\tt Map} files do define chain maps, i.e., that
the maps $m$ they specify do satisfy $dm=md$.

\end{enumerate}

\section{Machine processing of the data}
\label{sec:machine}

Modern computer languages are adept at processing text.   Nonetheless,
the raw data which is used to produce {\tt all.products} and {\tt
brackets.sym} is provided in the files {\tt s\_g/Map.aug} and {\tt
s\_g/brackets}, since this raw data may be easier to process by a
computer program.  Examples should suffice to make clear the
translation.

The entry which is reported as
\begin{verbatim}
  2    4  (  1    3     F2)  1_0
\end{verbatim}
in {\tt all.products} is derived from the line in {\tt 1\_0/Map.aug} which 
says
\begin{verbatim}
2 4 3
\end{verbatim}
as the rest of the data can be deduced from the map {\tt 1\_0}.

Similarly, the entries 
\begin{verbatim}
2_8 in < h4, 0, 1_0 >
2_1 in < h0, 1, 1_0 >
2_5 in < h3, 0, 1_0 >
\end{verbatim}
in {\tt 1\_0/brackets.sym} are derived from the lines in {\tt 1\_0/brackets} 
which say
\begin{verbatim}
2 8 16 0
2 1 1 1
2 5 8 0
\end{verbatim}
Here, the third entries in each line, $16$, $1$ and $8$, are the internal degrees
of the elements $h_4$, $h_0$ and $h_3$.

Finally, the entries in {\tt Sq0/Map.aug} have the form
\begin{verbatim}
2 1 0
2 3 1
2 5 2
\end{verbatim}
meaning that $Sq^0(2_0)$ contains $2_1$, that 
$Sq^0(2_1)$ contains $2_3$, and that 
$Sq^0(2_2)$ contains $2_5$.

\end{document}